\documentclass[preprint,3p,12pt]{elsarticle}

\usepackage{booktabs}
\usepackage{color}
\usepackage{amsmath}
\usepackage{amsfonts}
\usepackage{amssymb}
\usepackage{latexsym}
\usepackage{latexsym}
\usepackage{array}
\usepackage{amsthm}
\usepackage{graphics}
\usepackage{comment}
\usepackage{epsf,graphicx,subfigure}
\usepackage{epsfig}
\usepackage{epstopdf}
\usepackage{tikz}
\usetikzlibrary{arrows,positioning,shapes.geometric}

\usepackage{longtable}
\usepackage{rotating}

\usepackage{lscape}

\newtheorem{Theorem}{Theorem}
\newtheorem{Proposition}{Proposition}
\newtheorem{Lemma}{Lemma}

\newtheorem{Assumption}{Assumption}
\newtheorem{Remark}{Remark}
\newtheorem{Definition}{Definition}

\makeatletter
\def\ps@pprintTitle{%
  \let\@oddhead\@empty
  \let\@evenhead\@empty
  \let\@oddfoot\@empty
  \let\@evenfoot\@oddfoot
}
\makeatother

\begin{document}

\begin{frontmatter}

\title{Asymptotic Analysis of Higher-order Scattering Transform of Gaussian Processes}

\author[GRLiu]{Gi-Ren Liu}
\author[YCSheu]{Yuan-Chung Sheu}
\author[HTWu]{Hau-Tieng Wu\corref{mycorrespondingauthor}}
\cortext[mycorrespondingauthor]{Hau-Tieng Wu}
\ead{hauwu@math.duke.edu}

\address[GRLiu]{Department of Mathematics, National Chen-Kung University, Tainan, Taiwan}
\address[YCSheu]{Department of Applied Mathematics, National Yang Ming Chiao Tung University, Hsinchu, Taiwan}
\address[HTWu]{Department of Mathematics and Department of Statistical Science, Duke University, Durham, NC, USA}

\begin{abstract}
We analyze the scattering transform with the quadratic nonlinearity (STQN) of Gaussian processes without depth limitation. STQN is a nonlinear transform that involves a sequential interlacing convolution and nonlinear operators, which is motivated to model the deep convolutional neural network.
We prove that with a proper normalization, the output of STQN converges to a chi-square process with one degree of freedom in the finite dimensional distribution sense, and we provide a total variation distance control of this convergence at each time that converges to zero at an exponential rate.
To show these, we derive a recursive formula to represent the intricate nonlinearity of STQN by a linear combination of Wiener chaos, and then apply the Malliavin calculus and Stein's method to achieve the goal.
\end{abstract}


\end{frontmatter}


\section{Introduction}\label{sec:introduction}

The scattering transform (ST) is motivated by establishing a mathematical foundation of the convolutional neural network \cite{mallat2012group}, and it has been applied to various signals, for example, fetal heart rate \cite{chudavcek2013scattering}, brain waves \cite{liu2020diffuse,liu2020hospitals}, respiration \cite{wu2014assess}, marine bioacoustics \cite{balestriero2017linear}, and audio \cite{anden2011multiscale,li2019heart}. It provides a variety of representations for a given function $X$ through a sequential interlacing convolution and nonlinear operators:
\begin{equation*}
U[j_{1},j_{2},\ldots,j_M]X(t)=A\left(\cdots A\left(A\left(X\star\psi_{j_{1}}\right)\star\psi_{j_{2}}\right)\cdots \star \psi_{j_M}(t)\right),\ t\in\mathbb{R},
\end{equation*}
where {$M\in \mathbb{N}$ is the depth of ST,} $\{j_{1},j_{2},\ldots,j_{M}\}\subset \{\ldots,J-2,J-1,J\}$ is a set of scale parameters, $J\in \mathbb{Z}$ {determines the range of interest in the frequency domain},  $\{\psi_{j_{1}},\ldots,\psi_{j_{M}}\}$ is a family of wavelets generated from a selected mother wavelet $\psi$, and $A:\mathbb{C}\to\mathbb{C}$ is the chosen activation function. This construction is called the {\em scattering network}.
The associated $M$-th order ST coefficients are then computed through the pooling process
\begin{equation*}
S_{J}[j_{1},j_{2},\ldots,j_M]X(t) := U[j_{1},j_{2},\ldots,j_M]X\star \phi_{J}(t),\ t\in \mathbb{R},
\end{equation*}
where $\phi_{J}$ is a low-pass filter{, and it is usually chosen to be the father wavelet associated with $\psi$ and $J$ so that the Littlewood-Paley condition is satisfied. We mention that ST can also be defined with the Gabor transform \cite{czaja2019analysis,kavalerov20203} and others \cite{wiatowski2015deep,cheng2016deep}, but we focus on the wavelet transform considered in \cite{mallat2012group} in this paper}.
In most applications, only the first- and second-order ST coefficients are used because extracting higher-order information brings additional computational costs, except \cite{anden2011multiscale,anden2014deep}, which shows the potential benefit of considering higher order ST coefficients.
There are various choices for the activation $A$; for example, $A(\cdot)=|\cdot|$ \cite{mallat2012group,anden2014deep}, $A(\cdot) = |\cdot|^{2}$ \cite{balestriero2017linear,balan2018lipschitz}, and more general functions like Lipschtiz-continuous functions \cite{wiatowski2017mathematical,liu2020central}.
We mention that $A(\cdot) = |\cdot|^{2}$ is considered to speed up the computation of deeper scattering networks \cite{balestriero2017linear}, since the computation could be carried out without leaving the Fourier domain. When $A(\cdot)=|\cdot|^2$, we call the resulting transform the {\em ST with the quadratic nonlinearity} (STQN).

There have been several theoretical supports established for ST. If the mother wavelet and the low-pass filter $\phi_{J}$ satisfy the Littlewood-Paley condition and $A(\cdot) = |\cdot|$, the ST coefficients are approximately invariant to time shifts and stable to small time-warping deformation \cite{mallat2012group}. These properties have been generalized to Lipschitz-continuous activation functions \cite{liu2020central}. Similar time invariance and stability to deformation results are discussed in \cite{wiatowski2017mathematical} under a different setting.
The authors of \cite{balan2018lipschitz} prove that the STQN coefficients satisfy the Lipschitz continuous property; that is,
$\|S_{J}[j_{1},j_{2},\ldots,j_{M}]X-S_{J}[j_{1},j_{2},\ldots,j_{M}]Y\|_{L^{2}}$
can be controlled by $\|X-Y\|_{L^{2}}$ for any $X,Y\in L^{2}$, if
some conditions on the wavelets, the low-pass filter and the $L^{\infty}$ norm of inputs hold. The authors of \cite[Figure 4]{lostanlen2021one} showed that the decay of the STQN coefficients depend on the number of sinusoidal components contained in the signal and STQN has the ability to disentangle multiple factors of variability in the spectral envelope.
In addition to the above results in the deterministic setup, ST has also been studied from the probabilistic point of view. The authors of \cite{bruna2015intermittent} found that the second-order ST with $A(\cdot)=|\cdot|$ can be applied to characterize more detailed properties of random processes with stationary increments, by computing the expectation of the second-order ST of $X$ and observing its behavior with respect to the scale parameters $j_{1}$ and $j_{2}$. This result has been extended to general Lipschtiz-continuous activation functions \cite{liu2020central}.

Despite the above results, however, due to the nonlinearity of the activation function and its subsequent convolution, the behavior of higher order ST with a random process input is still open. Precisely, when $M>2$, since the input $U[j_{1},j_{2},\ldots,j_{M-1}]X$ to the $M$-th layer of the scattering network is neither a Gaussian process nor a subordinated Gaussian process even when $X$ is a stationary Gaussian process or a Gaussian process with stationary increments, it is challenging to study ST with random processes when $M>2$.
Motivated by its practical usefulness and theoretical challenges, in this work we focus on exploring STQN with depth $M\geq2$ from the probabilistic point of view, in which the inputs are the stationary Gaussian processes, or Gaussian processes with stationary increments.

To this end, several analysis techniques are needed. First, we derive a recursive formula to handle the intricate nonlinearity of STQN by representing $U[j_{1},j_{2},\ldots,j_{M}]X$ by a linear combination of Wiener chaos that comes from Wiener-It$\hat{\textup{o}}$ integrals of functions defined iteratively from the spectral density function of the Gaussian random process
Then, we apply the Stein's method \cite{nourdin2010invariance} and the product formula for Wiener-It$\hat{\textup{o}}$ integrals \cite{major1981lecture} to prove that
the total variation (TV) distance between $\sigma_{j_{M}}^{-2}U[j_{1},j_{2},\ldots,j_{M}]X(2^{j_{M}}t)$ and the square of the standard normal random variable is $O(2^{-j_{M}/2})$ for each fixed $t\in \mathbb{R}$, where $\sigma_{j_{M}}^{2} = \mathbb{E}\left[U[j_{1},j_{2},\ldots,j_{M}]X(\cdot)\right]$. Moreover, we prove that $2^{j_{M}}U[j_{1},j_{2},\ldots,j_{M}]X(2^{j_{M}}\cdot)$ converges to a weakly dependent chi{-}square process with one degree of freedom in the finite dimensional distribution sense as $j_{M}\rightarrow\infty$.
We shall mention that the proof in this work is inspired by \cite{nourdin2009stein}, where the Malliavin calculus and Stein's method are combined to derive explicit bounds in the Gaussian and Gamma approximations of random variables in a fixed Wiener chaos of a general Gaussian process. Although the work \cite[Proposition 3.7]{nourdin2009stein} provides an error bound for the Gaussian approximation of a finite sum of multiple Wiener-It$\hat{\textup{o}}$ integrals, it cannot be directly applied to our problem. The main obstacle comes from the intricate covariance structure between different orders of Wiener chaos, which need new controls to ensure that the approximation error converges to zero as $j_{M}\rightarrow\infty$.
%

The rest of the paper is organized as follows.
In Section \ref{sec:preliminary}, we {summarize necessary material for ST and} present some preliminaries about the Wiener-It$\hat{\textup{o}}$ integrals
and the Malliavin calculus. In Section \ref{sec:mainresult}, we state our main
results, including Theorems \ref{thm:conv_rate_4rd},
\ref{thm:equiv_F}, and
\ref{thm:equiv_F_plus}.
The proofs of our main results and some technical lemmas are given
in Section \ref{sec:proof}. Table \ref{List_symbols} below contains a list of frequently used symbols and abbreviations.

\begin{table}\label{List_symbols}
\begin{center}\begin{tabular}{p{2.7cm}p{13cm}}
\toprule
{\bf Symbol} & {\bf Description}
\\
\midrule
$\star$  & Convolution operator\\
$M$      & Depth of the scattering network\\
$\pm\lambda_{1:p}$    & Abbreviation for $\pm(\lambda_{1},\lambda_{2},\ldots,\lambda_{p})\in \mathbb{R}^{p}$\\
$\lambda^{+}_{1:p}$, $d\lambda_{1:p}$,      & Abbreviation for $\lambda_{1}+\lambda_{2}+\cdots+\lambda_{p}$ and $d\lambda_{1}\ d\lambda_{2}\cdots d\lambda_{p}$\\
$j_{1:m}$ & Vector of scale parameters for the wavelet transform in the first $m$ layers\\
$\psi$, $\hat{\psi}$, $\psi_{j}(\cdot)$ & Mother wavelet, its Fourier transform, and its scaled version $2^{-j}\psi(2^{-j}\cdot)$\\
$X$       & Stationary Gaussian process\\
$U[j]X$   & Modulus-squared wavelet transform of  $X$ defined by  $U[j]X = |\psi_{j}\star X|^{2}$\\
$U[j_{1:m}]$ & The $m$-th order scattering transform defined by $U[j_{m}]U[j_{m-1}]\cdots U[j_{1}]$\\
$T[j_{1:m+1}]X$ &  $\left(U[j_{1:m}]X\right)\star \psi_{j_{m+1}}$\\
$\alpha$  &  Vanishing-moment parameter for $\psi$ defined through $\hat{\psi}(\cdot) = C_{\hat{\psi}}(\cdot)|\cdot|^{\alpha}$\\
$f_{Z}$   & Spectral density of the covariance function of a random process $Z$\\
$\beta$       & Singularity parameter for $f_{X}$ defined through $f_{X}(\cdot) = C_{X}(\cdot)|\cdot|^{\beta-1}$\\
$\overline{H}$ & Hilber space $\{f\in L^{2}(\mathbb{R})\mid f(-\lambda) = f(\lambda)\ \textup{for all}\ \lambda\in \mathbb{R}\}$\\
$\overline{H}^{\otimes p}$ ($\overline{H}^{\odot p}$) & The $p$-th (symmetric) tensor product of the Hilbert space $\overline{H}$\\
$W$  & Complex-valued Gaussian random measure\\
$I_{p}$       & $p$-fold Wiener-It$\hat{\textup{o}}$ integrals for $p\in \mathbb{N}$ ($I_{0}$: the identity map)\\
$\widetilde{f}$       & Canonical symmetrization of function $f$\\
$f\otimes_{r} g$, $f\widetilde{\otimes}_{r} g$    & The $r$-th contraction of $f$ and $g$, and its canonical symmetrization\\
$L^2(\Omega,\overline{H}^{\otimes p})$  & Random functions $f$ in the Hilbert space $\overline{H}^{\otimes p}$ with $\mathbb{E}[\|f\|^{2}_{L^{2}(\mathbb{R}^{p})}]<\infty$\\
$M_{r}$     & $2^{M-1}-2r$ for $r=0,1,\ldots,2^{M-2}-1$\\
$\{{}^{M}F^{(M_{r})}_{t}\}_{r=0}^{2^{M-2}-1}$  & Integrands in the Wiener-It$\hat{\textup{o}}$ decomposition
of $T[j_{1:M}]X$\\
$\{{}^{M}\widetilde{G}^{(2^{M}-2\ell)}_{t}\}_{\ell=0}^{2^{M-1}}$  & Integrands in the Wiener-It$\hat{\textup{o}}$ decomposition
of $U[j_{1:M}]X$\\
$a_{n} = O(b_{n})$  & For sequences $\{a_{n}\}$ and $\{b_{n}\}$, there exists
a constant $C>0$ such that $|a_{n}|\leq C |b_{n}|$ when $n$ is sufficiently large.\\
\bottomrule
\end{tabular}
\end{center}
\caption{List of frequently used symbols and abbreviations}
\end{table}

\vskip 20 pt

\section{Preliminaries}\label{sec:preliminary}

\subsection{Scattering transform}
Let $\psi\in L^{1}(\mathbb{R})\cap L^{2}(\mathbb{R})$ be the mother wavelet, which satisfies $\int_{\mathbb{R}} \psi(t) dt = 0$.
A family of real-valued functions $\{\psi_{j}(t)\mid j\in \mathbb{Z}, t \in \mathbb{R}\}$
is called a wavelet family, {where} 
\begin{equation*}
\psi_{j}(t) = \frac{1}{2^{j}}\psi\left(\frac{t}{2^{j}}\right)\,.
\end{equation*}
Denote the Fourier transform of $\psi$ by $\hat{\psi}$, which is defined as
\begin{equation*}
\hat{\psi}(\lambda) = \int_{\mathbb{R}}e^{-i\lambda t}\psi(t) dt.
\end{equation*}
Because $\psi\in L^{1}(\mathbb{R})\cap L^{2}(\mathbb{R})$ and $\int_{\mathbb{R}} \psi(t) dt = 0$,
$\hat{\psi}\in {C^0(\mathbb{R})\cap L^2(\mathbb{R})}$ with $\hat{\psi}(0) = 0$. {We} make the following assumption.

\begin{Assumption}\label{Assumption:1:wavelet}
{We} assume that $\hat{\psi}\in L^{1}(\mathbb{R})$ and there exists a bounded and continuous {complex-valued} function $C_{\hat{\psi}}$, which is positive at the origin (i.e., $C_{\hat{\psi}}(0)>0$),
such that $\hat{\psi}$ can be expressed as
\begin{equation*}
\hat{\psi}(\lambda) = C_{\hat{\psi}}(\lambda)|\lambda|^{\alpha}
\end{equation*}
for some $\alpha>0$.
\end{Assumption}

\begin{Remark}
{Consider} the vanishing moment $N$ of $\psi$ {defined as}
\begin{equation*}
N = \max\left\{n\in \mathbb{N} \mid \int_{\mathbb{R}}t^{\ell}\psi(t)dt=0\ \textup{for}\ \ell = 0,1,\ldots,n-1,\ \textup{and}\  \int_{\mathbb{R}}t^{n}\psi(t)dt\neq0\right\}.
\end{equation*}
{The parameter $\alpha$ in Assumption \ref{Assumption:1:wavelet}
satisfies $\alpha\geq N$ by the Taylor expansion.}
Assumption \ref{Assumption:1:wavelet} with $\alpha\geq 1$ holds for commonly used wavelets, including the real part of complex Morlet wavelet ($\alpha=1$), the Mexican hat wavelet ($\alpha=2$), and the Daubechies-$K$ wavelet ($\alpha=K/2$), where $K = 2,4,\ldots,20.$
{We shall mention that for those} wavelet functions whose Fourier transform $\hat{\psi}$ vanishes in a neighborhood of the origin, e.g., the Mayer wavelet,
all results in this paper hold.
\end{Remark}

Given a locally bounded function $X: \mathbb{R}\rightarrow \mathbb{R}$, the wavelet transform of $X$ is defined as
\begin{equation}\label{df:Wavelet}
      X\star \psi_{j} (t) = \int_{\mathbb{R}}X(s)\psi_{j}(t-s)ds,\ j\in \mathbb{Z},\ t\in \mathbb{R}.
\end{equation}
The first-order {STQN} of $X$ is obtained by composing
the wavelet transform of $X$ with the nonlinear operator $|\cdot|^{2}$; {that is,}
\begin{equation*}
U[j_{1}]X(t) = |X\star \psi_{j_{1}}(t)|^{2},\ j_{1}\in \mathbb{Z}.
\end{equation*}
Fix $M\geq 1$ to be the order (or depth) of ST. The $M$-th order STQN of $X$ is defined iteratively as
$$
U[j_{1},j_{2},\ldots,j_{M}]X(t) := U[j_{M}]U[j_{M-1}]\cdots U[j_{2}]U[j_{1}]X(t)=: U[j_{1:M}]X(t)\,,
$$
where $j_{1},j_{2},\ldots,j_{M}\in \mathbb{Z}$ and $j_{1:M}:=[j_{1},\ldots,j_{M}]\in \mathbb{Z}^{M}$ is introduced to simplify the heavy notation for the upcoming proof.
Denote the wavelet transform of the output of the $M$-th order STQN of $X$ by
\begin{equation}\label{def:T}
T[j_{1:M+1}]X= U[j_{1:M}]X\star\psi_{j_{M+1}}\,.
\end{equation}
We define $T[j_{1}]X= X\star\psi_{j_{1}}$. By definition, the operators $U$ and $T$ have the relationship
\begin{equation}\label{U=squareT}
U[j_{1:M}]X = \left|T[j_{1:M}]X\right|^{2}.
\end{equation}


\subsection{Stationary Gaussian processes}
In this section, we prepare materials about stationary Gaussian random processes.
Let $\overline{H}=\{f\in L^{2}(\mathbb{R})\mid f(-\lambda)=\overline{f(\lambda)}\ \textup{for all}\ \lambda\in \mathbb{R}\}$ be a complex Hilbert space
with {the} inner product
$
\langle f,g\rangle = \int_{\mathbb{R}}f(\lambda)\overline{g(\lambda)}d\lambda,
$
and let $W$ be a complex-valued Gaussian random measure on $\mathbb{R}$
satisfying
\begin{align}\label{ortho}
W(\Delta_{1})=\overline{W(-\Delta_{1})},\ \
\mathbb{E}[W(\Delta_{1})]=0,\
\textup{and}\ \ \mathbb{E}\left[W(\Delta_{1})
\overline{W(\Delta_{2})}\right]=\textup{Leb}(\Delta_{1}\cap\Delta_{2})
\end{align}
for any
 $\Delta_{1},\Delta_{2}\in
\mathcal{B}(\mathbb{R})$, where Leb is the Lebesgue measure on $\mathbb{R}$ {and $\mathcal{B}(\mathbb{R})$ is the
Borel {$\sigma$-}algebra on $\mathbb{R}$}.
For every $f\in\overline{H}$, we define $W(f) = \int_{\mathbb{R}}f(\lambda)W(d\lambda).$
{Due to} (\ref{ortho}),
$W(f)\in \mathbb{R}$ and $\mathbb{E}[W(f)W(g)] = \langle f,g\rangle.$
Hence, $W=\{W(f)\mid f\in \overline{H}\}$ can be viewed as an centered isonormal Gaussian process over $\overline{H}$
{\cite[Section 2.1]{nourdin2012normal}}.
Let $(\Omega, \mathcal{F}_{W}, {\mathcal P})$ be a
probability space, where the $\sigma$-algebra $\mathcal{F}_{W}$ is generated by $W$.

Let $X$ be a mean-square continuous and
stationary real Gaussian random process with the constant mean $\mathbb E[X(s)]=\mu\in \mathbb{R}$ and the
covariance function $R_{X}$:
\begin{equation*}
R_{X}(t_{1}-t_{2}) =\mathbb E[X(t_{1})X(t_{2})]-\mathbb E[X(t_{1})]\mathbb E[X(t_{2})],\ t_{1},t_{2}\in \mathbb{R}.
\end{equation*}
{Note that $R_X$ is continuous by the assumption of $X$. Since} $\int_{\mathbb{R}}\psi(s) ds = 0 $,
$X\star\psi_{j}(t)
= [X-\mu]\star\psi_{j}(t)\,,
$
{and hence} we can assume that $\mathbb E[X(s)]=0$ without loss of generality.
By the Bochner-Khinchin theorem,
there exists a unique nonnegative measure $F_{X}:\mathcal{B}(\mathbb{R})\rightarrow [0,\,\infty)$ such that $F_{X}(\Delta) = F_{X}(-\Delta)$ for any $\Delta\in\mathcal{B}(\mathbb{R})$
and
\begin{equation}\label{eq:BK}
R_{X}(t) = \int_{\mathbb{R}}e^{i\lambda t}F_{X}(d\lambda),\ t\in \mathbb{R}.
\end{equation}
The measure $F_{X}$ is called the spectral measure of the covariance function $R_{X}$.

\begin{Assumption}\label{Assumption:2:spectral}
The spectral measure $F_{X}$ has the density
$f_{X}$ and
\begin{equation*}
f_{X}(\lambda) = \frac{C_{X}(\lambda)}{|\lambda|^{1-\beta}},
\end{equation*}
where $\beta\in(0,1)$ is the {\em Hurst index} of long-range {dependence} and $C_{X}$ is a bounded and continuous function from $\mathbb{R}$ to $[0,\infty)$
such that $C_{X}(\lambda)$ {decays} faster than $|\lambda|^{-\beta-\varepsilon}$ for some $\varepsilon>0$ {when} $|\lambda|\rightarrow\infty$.
\end{Assumption}

If the function $C_{X}$ in Assumption \ref{Assumption:2:spectral} satisfies $C_{X}(0)>0$, then $X$ is a long-range dependent process \cite{doukhan2002theory,pipiras2017long} because
$f_{X}$ has a singularity at {$0$}. Note that the spectral density function $f_{X}$ is even and nonnegative with $\int_{\mathbb{R}} f_{X}(\lambda) d\lambda = R_{X}(0)$.

Under Assumption \ref{Assumption:2:spectral}, (\ref{eq:BK}) can be rewritten as
\begin{equation*}
R_{X}(t) = \int_{\mathbb{R}}e^{i\lambda t}f_{X}(\lambda)d\lambda,\ t\in \mathbb{R}.
\end{equation*}
By the Karhunen theorem, the Gaussian process $X$ has the
representation
\begin{align}\label{sample path represent}
X(t)=\int_{\mathbb{R}}
e^{i\lambda t}\sqrt{f_{X}(\lambda)}W(d\lambda),\ t\in\mathbb{R}.
\end{align}
%
%
By substituting (\ref{sample path represent}) into (\ref{df:Wavelet}),
the continuous wavelet transform of $X$ can be represented as a Wiener-It$\hat{\textup{o}}$ integral
\begin{align}\notag
X\star\psi_{j}(t) =& \int_{\mathbb{R}}\left[\int_{\mathbb{R}}
e^{i\lambda s}\sqrt{f_{X}(\lambda)}W(d\lambda)\right]\psi_{j}(t-s)ds
\\=&\int_{\mathbb{R}}\left[\int_{\mathbb{R}}
e^{i\lambda s}\psi(t-s)ds\right]\sqrt{f_{X}(\lambda)}W(d\lambda)
\label{lemma:scatter1}=\int_{\mathbb{R}}
e^{i\lambda t}\hat{\psi}_{j}(\lambda)\sqrt{f_{X}(\lambda)}W(d\lambda),
\end{align}
where the last equality follows from the stochastic Fubini theorem \cite[Theorem 2.1]{pipiras2010regularization} (see also \ref{appendix:fubini} for more details).
From (\ref{lemma:scatter1}), we know that $X\star\psi_{j}$ has the spectral density
$$
f_{X\star \psi_{j}}(\lambda) =  f_{X}(\lambda)|\hat{\psi}_{j}(\lambda)|^{2},\ \lambda\in \mathbb{R},
$$
and
\begin{equation}\label{U1}
U[j_{1}]X(t) = \left(\int_{\mathbb{R}}
e^{i\lambda t}\hat{\psi}_{j_{1}}(\lambda)\sqrt{f_{X}(\lambda)}W(d\lambda)\right)^{2}.
\end{equation}

\begin{Remark}
It is known since the works of It$\hat{\textup{o}}$ \cite{ito1954stationary} and Yaglom \cite{yaglom1957some} that if a
Gaussian random process $X$ is increment-stationary in the sense that for any $s$, $s'$, $t$, $t'$ and
$h\in \mathbb{R}$,
\begin{equation*}
\mathbb{E}\left[\left(X(t+h)- X(s+h)\right)\left(X(t'+h)- X(s'+h)\right)\right] = \mathbb{E}\left[\left(X(t)-X(s)\right)\left(X(t') - X(s') \right)\right]
\end{equation*}
and if
the covariance between the increments over time intervals $I$ and $I'$ tends to zero when
$\min\{|a-b|\mid a\in I,\ b\in I^{'}\}\rightarrow\infty$,
then $X$ admits a spectral representation
\begin{align}\label{spectral_stat_increment}
X(t) = c+\int_{\mathbb{R}}\left(e^{it\lambda}-1\right)Z(d\lambda),\ t\in \mathbb{R},
\end{align}
where $c\in \mathbb{R}$ and  $Z(d\lambda)$ is a
complex-valued Gaussian random measure on $\mathbb{R}$ with spectral measure $F_{X}$.
If the spectral measure $F_{X}$ has a density $f_{X}$, then (\ref{spectral_stat_increment})
can be rewritten as
\begin{align*}
X(t) = c+\int_{\mathbb{R}}\left(e^{it\lambda}-1\right)\sqrt{f_{X}(\lambda)}W(d\lambda),\ t\in \mathbb{R}.
\end{align*}
By the mean-zero property of the wavelet,
\begin{align}\notag
X\star \psi_{j}(t) =&\int_{\mathbb{R}}\left[ \int_{\mathbb{R}}\left(e^{it'\lambda}-1\right)\sqrt{f_{X}(\lambda)}W(d\lambda) \right]\psi_{j}(t-t')dt'
\\\label{wavelet_stat_increment}=& \int_{\mathbb{R}}e^{it\lambda}\hat{\psi}_{j}(\lambda)\sqrt{f_{X}(\lambda)}W(d\lambda).
\end{align}
From the consistency between (\ref{lemma:scatter1}) and (\ref{wavelet_stat_increment}),
we know that the all results below hold for Gaussian random processes with stationary increments, including the fractional Brownian motions, when the associated conditions are fulfilled.
\end{Remark}

\subsection{Wiener-It$\hat{o}$ integrals and Winer chaos expansion}
In this section, we describe how to express the STQN of stationary Gaussian random processes in terms of the Wiener chaos expansion that comes from Wiener-It$\hat{\textup{o}}$ integrals of functions defined iteratively from the spectral density function of the Gaussian random process.
Given an integer $p\geq2$, the $p$-th tensor product of the Hilbert space $\overline{H}$ is denoted by $\overline{H}^{\otimes p}$.
{According to \cite[p. 27]{major1981lecture},} $f\in \overline{H}^{\otimes p}$ if and only if $f=f(\lambda_{1},\ldots,\lambda_{p})$,
where $(\lambda_{1},\ldots,\lambda_{p})\in \mathbb{R}^{p}$, is a complex valued function
satisfying $f(-\lambda_{1},\ldots,-\lambda_{p})=\overline{f(\lambda_{1},\ldots,\lambda_{p})}$
and
$$\|f\|^{2}_{\overline{H}^{\otimes p}}  =\int_{\mathbb{R}^{p}}|f(\lambda_{1:p})|^{2}d\lambda_{1:p} <\infty\,,
$$
where $\lambda_{p_1:p_2}$ means $(\lambda_{p_1},\,\lambda_{p_1+1}\ldots,\lambda_{p_2})\in \mathbb{R}^{p_2-p_1+1}$ and $d\lambda_{p_1:p_2}$ means $d\lambda_{p_1}d\lambda_{p_1+1}\cdots d\lambda_{p_2}$ when $p_2>p_1$.
For any $f,g\in \overline{H}^{\otimes p}$, their inner
product is defined as
\begin{align}\notag
\langle f,g\rangle_{\overline{H}^{\otimes p}} = \int_{\mathbb{R}^{p}}f(\lambda_{1:p})\overline{g(\lambda_{1:p})}
d\lambda_{1:p}
=\int_{\mathbb{R}^{p}}f(\lambda_{1:p})g(-\lambda_{1:p})
d\lambda_{1:p}\,.
\end{align}
The subscripts of $\|f\|^{2}_{\overline{H}^{\otimes p}}$ and $\langle f,g\rangle_{\overline{H}^{\otimes p}}$
will be ignored as {they} can be inferred from the number of variables of $f$ and $g$.
The $p$-th symmetric tensor product of $\overline{H}$
is denoted by $\overline{H}^{\odot p}$, which contains those functions $f\in \overline{H}^{\otimes p}$
satisfying
$f(\lambda_{\pi(1)},\ldots,\lambda_{\pi(p)}) = f(\lambda_{1},\ldots,\lambda_{p})$, where $\pi$ {is any permutation} of the set $\{1,2,\ldots,p\}$.
Because $\overline{H}^{\odot p}$ {is a Hilbert subspace of} $\overline{H}^{\otimes p}$,
the norm and inner product in $\overline{H}^{\odot p}$ are defined in the same way as in $\overline{H}^{\otimes p}$.
For any $f\in \overline{H}^{\otimes p}$, where $p\in \mathbb{N}$,
the {\em $p$-fold Wiener-It$\hat{\textup{o}}$ integrals} of $f$ with respect to the random measure $W$ is defined by
\begin{align*}
I_{p}(f) = \int_{\mathbb{R}^{p}}^{'}f(\lambda_{1:p})W(d\lambda_{1})\cdots W(d\lambda_{p}),
\end{align*}
where
$\int^{'}$ means that the integral excludes the diagonal hyperplanes
$\lambda_{k}=\mp \lambda_{k^{'}}$ for $k, k^{'}\in\{1,\ldots,p\}$ and  $k\neq k^{'}$.
Note that $I_{p}(f) = I_{p}(\tilde{f})$, where
$\tilde{f}\in \overline{H}^{\odot p}$ is the canonical symmetrization of $f$ defined as
\begin{align}\label{def:canonical_symmetrization}
\tilde{f}(\lambda_{1:p}) \frac{1}{p!}\underset{\pi}{\sum}f(\lambda_{\pi(1)},\ldots,\lambda_{\pi(p)}),
\end{align}
where the sum runs over all permutations $\pi$ of $\{1,\ldots,p\}$.
By default, $I_{0}$ is the identity map, i.e., $I_{0}(f) = f$ for any function $f$.

For $p\in \mathbb{N}\cup\{0\}$, let $H_p$ be the Hermite polynomial of degree $p$, which is defined
by
\[
{h}_{p}(y)=(-1)^{p}e^{\frac{y^{2}}{2}}\frac{d^{p}}{dy^{p}}e^{-\frac{y^{2}}{2}},\ y\in \mathbb{R}.
\]
We write $\mathcal{H}_{p}$ to denote the closed linear subspace of $L^{2}(\Omega,\mathcal{F}_{W},\mathcal{P})$
generated by the random variables of type ${h}_{p}(W(g))$, $g\in \overline{H}$, $\|g\|=1$.
The space $\mathcal{H}_{p}$ is called the $p$-th Wiener chaos of $W$.
According to \cite{major1981lecture,nualart2006malliavin,nourdin2009stein,nourdin2012normal}, $I_{p}(\tilde{f})\in \mathcal{H}_{p}$ for any $f\in \overline{H}^{\otimes p}$.
The following {product formula} allows us to express output of STQN in terms of a linear combination of Wiener chaos, which is the first  tool used in this paper.

\begin{Lemma}[Product Formula \cite{major1981lecture,nourdin2012normal}]\label{lemma:itoformula}
Let $p,q\geq 1$. If $f\in \overline{H}^{\odot p}$ and $g\in \overline{H}^{\odot q}$, then
\begin{align*}
I_{p}(f)I_{q}(g) = \overset{p\wedge q}{\underset{r=0}{\sum}}r!\binom{p}{r}\binom{q}{r}I_{p+q-2r}\left(f\otimes_{r} g\right),
\end{align*}
where $f\otimes_{r} g$ is
the $r$th contraction of $f$ and $g$ defined as
\begin{align*}
f\otimes_{r} g(\lambda_{1:p+q-2r})
= \int_{\mathbb{R}^{r}} f(\tau_{1:r},\lambda_{1:p-r})
g(-\tau_{1:r},\lambda_{p-r+1:p+q-2r})d\tau_{1:r}
%
\end{align*}
for $r=1,2,\ldots,p\wedge q$.
{When $r=0$, set} $f\otimes_{0} g := f\otimes g$.

\end{Lemma}
Note that $f\otimes_{r} g\in\overline{H}^{\otimes p+q-2r}$, but it may not be symmetric.
Let $f\tilde{\otimes}_{r} g$ be the canonical symmetrization of $f\otimes_{r} g$.
{
For any $f\in \overline{H}^{\odot p}$,
by noting that the expectation of any Wiener-It$\hat{\textup{o}}$ integral equals to zero,
Lemma \ref{lemma:itoformula} implies
\begin{align*}
\mathbb{E}\left[|I_{p}(f)|^{2}\right]
=
\mathbb{E}\left[
\overset{p}{\underset{r=0}{\sum}}r!\binom{p}{r}^{2}I_{2p-2r}\left(f\otimes_{r} f\right)\right]
= p!f\otimes_{p}f = p!\|f\|_{L^{2}(\mathbb{R}^{p})}^{2},
\end{align*}
which is the so-called
isometry property of the multiple Wiener-It$\hat{\textup{o}}$ integrals.}

{We now express the ST coefficients using the product formula.} By applying Lemma \ref{lemma:itoformula} to (\ref{U1}),
\begin{equation}\label{U1:wienerchaos}
U[j_{1}]X(t) = \int^{'}_{\mathbb{R}^{2}}
e^{i(\lambda_{1}+\lambda_{2}) t}\hat{\psi}_{j_{1}}(\lambda_{1})\sqrt{f_{X}(\lambda_{1})}\hat{\psi}_{j_{1}}(\lambda_{2})\sqrt{f_{X}(\lambda_{2})}W(d\lambda_{1})W(d\lambda_{2})
+\|f_{X\star \psi_{j_{1}}}\|_{1}.
\end{equation}
Denote
\begin{align}\label{def:1F2}
\left\{
\begin{array}{l}
{}^{1}\widetilde{G}^{(0)}_{t} = \|f_{X\star \psi_{j_{1}}}\|_{1},\\
{}^{1}\widetilde{G}^{(2)}_{t}(\lambda_{1},\lambda_{2}) = e^{i(\lambda_{1}+\lambda_{2}) t}\hat{\psi}_{j_{1}}(\lambda_{1})\sqrt{f_{X}(\lambda_{1})}\hat{\psi}_{j_{1}}(\lambda_{2})\sqrt{f_{X}(\lambda_{2})}.
\end{array}
\right.
\end{align}
The first-order STQN of $X$ can be rewritten as
\begin{equation*}
U[j_{1}]X(t) = I_{2}\left({}^{1}\widetilde{G}^{(2)}_{t}\right)+{}^{1}\widetilde{G}^{(0)}_{t}.
\end{equation*}
By the property $\int_{\mathbb{R}}\psi_{j_{2}}(t-t^{'})dt^{'} = 0$ and (\ref{U1:wienerchaos}),
\begin{align}\label{T2_wienerchaos_half}
T[j_{1:2}]X(t) =& \int_{\mathbb{R}}I_{2}\left({}^{1}\widetilde{G}^{(2)}_{t^{'}}\right)
 \psi_{j_{2}}(t-t^{'})dt^{'}.
\end{align}
Because (\ref{def:1F2}) shows that
\begin{align}\label{time_relation_1G}
{}^{1}\widetilde{G}^{(2)}_{t}(\lambda_{1},\lambda_{2}) = {}^{1}\widetilde{G}^{(2)}_{0}(\lambda_{1},\lambda_{2})\ e^{i(\lambda_{1}+\lambda_{2})t},
\end{align}
\begin{equation}\label{material_T2}
{}^{2}F^{(2)}_{t}(\lambda_{1},\lambda_{2}) :=
\int_{\mathbb{R}}{}^{1}\widetilde{G}^{(2)}_{t^{'}}(\lambda_{1},\lambda_{2})\psi_{j_{2}}(t-t^{'})dt^{'}
={}^{1}\widetilde{G}^{(2)}_{t}(\lambda_{1},\lambda_{2})
\hat{\psi}_{j_{2}}\left(\lambda_{1}+\lambda_{2}\right),
\end{equation}
which depends on $t$ in the same way.
More concretely,
$$
{}^{2}F^{(2)}_{t}(\lambda_{1},\lambda_{2})=e^{i(\lambda_{1}+\lambda_{2}) t}\hat{\psi}_{j_{1}}(\lambda_{1})\sqrt{f_{X}(\lambda_{1})}\hat{\psi}_{j_{1}}(\lambda_{2})\sqrt{f_{X}(\lambda_{2})}\hat{\psi}_{j_{2}}\left(\lambda_{1}+\lambda_{2}\right).
$$
By substituting (\ref{material_T2}) into (\ref{T2_wienerchaos_half}), we get
\begin{align}\label{T2_wienerchaos}
T[j_{1:2}]X(t) =& I_{2}\left({}^{2}F^{(2)}_{t}\right).
\end{align}
By the product formula {again} with $p=q=2$ and $f=g={}^{2}F^{(2)}_{t}$,
\begin{align}\label{U2_wienerchaos}
U[j_{1:2}]X(t)=\left[T[j_{1:2}]X(t)\right]^{2} = \left[I_{2}\left({}^{2}F^{(2)}_{t}\right)\right]^{2}
=
\overset{2}{\underset{r=0}{\sum}}I_{4-2r}\left({}^{2}\widetilde{G}^{(4-2r)}_{t}\right),
\end{align}
where
\begin{equation}\label{2G2F}
{}^{2}\widetilde{G}^{(4-2r)}_{t} = r!\binom{2}{r}\binom{2}{r}{}^{2}F^{(2)}_{t}\widetilde{\otimes}_{r} {}^{2}F^{(2)}_{t}\in \overline{H}^{\odot 4-2r}
\end{equation}
for each $t\in \mathbb{R}$. From (\ref{time_relation_1G}) and (\ref{material_T2}),
we have
\begin{align}\label{time_relation_2G}
{}^{2}\widetilde{G}^{(4-2r)}_{t}(\lambda_{1:4-2r}) = e^{i\lambda^{+}_{1:4-2r}t}\ {}^{2}\widetilde{G}^{(4-2r)}_{0}(\lambda_{1:4-2r})
\end{align}
for $r=0,1$, where $\lambda^+_{p_1:p_2}$ means $\lambda_{p_1}+\lambda_{p_1+1}+\ldots+\lambda_{p_2}\in \mathbb{R}$ for $p_2>p_1\geq 1$.
%
For general cases, we have the following result.
\begin{Lemma}\label{lemma:recursiveTU}
For every integer $M\geq2$ and stationary Gaussian process $X$, the random processes $T[j_{1:M}]X$ and $U[j_{1:M}]X$ arising from {STQN} of $X$
have the Wiener-It$\hat{o}$ decomposition as follows
\begin{align}\label{TM_wienerchaos}
T[j_{1:M}]X(t)
=\overset{2^{M-2}-1}{\underset{r=0}{\sum}}I_{2^{M-1}-2r}\left({}^{M}\!F^{(2^{M-1}-2r)}_{t}
\right)
\end{align}
and
\begin{align}\label{UM_wienerchaos}
U[j_{1:M}]X(t)
=\overset{2^{M-1}}{\underset{\ell=0}{\sum}}I_{2^{M}-2\ell}\left({}^{M}\widetilde{G}^{(2^{M}-2\ell)}_{t}
\right),
\end{align}
where the integrands, ${}^{M}F^{(2^{M-1}-2r)}_{t}$ and ${}^{M}\widetilde{G}^{(2^{M}-2\ell)}_{t}$, are defined below using the following recursive formula starting with (\ref{def:1F2}).
For $\ell = 0,1,\ldots,2^{M-1}$,
\begin{align}
{}^{M}\widetilde{G}^{(2^{M}-2\ell)}_{t}(\lambda_{1:2^{M}-2\ell})
:=&\overset{2^{M-2}-1}{\underset{r,r^{'}=0}{\sum}}
(\ell-r-r')!\binom{2^{M-1}-2r}{\ell-r-r'}\binom{2^{M-1}-2r'}{\ell-r-r'}\nonumber\\
&\times {}^{M}F^{(2^{M-1}-2r)}_{t}\widetilde{\otimes}_{\ell-r-r'} {}^{M}F^{(2^{M-1}-2r')}_{t}(\lambda_{1:2^{M}-2\ell}),\label{GM_wienerchaos}
\end{align}
where for $r = 0,1,\ldots,2^{M-2}-1$,
\begin{equation}\label{FM_wienerchaos}
{}^{M}F^{(2^{M-1}-2r)}_{t}(\lambda_{1:2^{M-1}-2r})
:=
{}^{M-1}\widetilde{G}^{(2^{M-1}-2r)}_{t}(\lambda_{1:2^{M-1}-2r})
\hat{\psi}_{j_{M}}(\lambda^{+}_{1:2^{M-1}-2r}).
\end{equation}
\end{Lemma}
For the binomial coefficients in (\ref{GM_wienerchaos}), we use the following rule: for any $A,B,C\in \mathbb{Z}$,
$$\binom{A}{C}\binom{B}{C}=0\ \textup{if}\
C<0
\ \textup{or}\
\min\{A,B\} <C.$$
The proof of Lemma \ref{lemma:recursiveTU} is based on the mathematical induction method, which can be found in Section \ref{sec:proof:lemma:recursiveTU}.
For the functions generated by the recursive formula (\ref{GM_wienerchaos}) and (\ref{FM_wienerchaos}) with the initial term (\ref{def:1F2}),
we have the following observation, which will be used in the proof of Lemma \ref{lemma:recursiveTU} and Proposition \ref{lemma:estimate_norm_fpqr}. The first one is
\begin{equation}\label{time_relation_MG}
{}^{M}\widetilde{G}^{(2^{M}-2\ell)}_{t}(\lambda_{1:2^{M}-2\ell}) = \textup{exp}\left(it \lambda^{+}_{1:2^{M}-2\ell}\right) \ {}^{M}\widetilde{G}^{(2^{M}-2\ell)}_{0}(\lambda_{1:2^{M}-2\ell})
\end{equation}
for $\ell = 0,1,\ldots,2^{M-1}-1$
and the second one is
\begin{equation}\label{time_relation_MF}
{}^{M}F^{(2^{M-1}-2r)}_{t}(\lambda_{1:2^{M-1}-2r}) = \textup{exp}\left(it \lambda^{+}_{1:2^{M-1}-2r}\right)\  {}^{M}F^{(2^{M-1}-2r)}_{0}(\lambda_{1:2^{M-1}-2r})
\end{equation}
for $r = 0,1,\ldots,2^{M-2}-1$, which are generalization of, for instance, (\ref{material_T2}) and (\ref{time_relation_2G}).

\tikzstyle{line} = [draw, -latex']
\tikzstyle{arrow} = [thick,->,>=stealth]

\begin{figure}[h!]  
\centering
   \begin{tikzpicture}[>=latex']
        \tikzset{
        block/.style= {draw, rectangle, align=center,minimum width=1.1cm,minimum height=.10cm,line width=0.05mm},
        rblock/.style={draw, shape=diamond,rounded corners=1.5em,align=center,minimum width=3cm,minimum height=.10cm},
        }

        \node []  (1G) {${}^{1}\widetilde{G}$};
        \node [, right =1.0cm of 1G]  (2F) {${}^{2}F$};
        \node [, right =1.0cm of 2F]  (2G) {${}^{2}\widetilde{G}$};
        \node [, right =1.0cm of 2G]  (3F) {${}^{3}F$};
        \node [, right =1.0cm of 3F]  (3G) {${}^{3}\widetilde{G}$};
        \node [, right =1.0cm of 3G]  (4F) {${}^{4}F$};
        \node [, right =1.0cm of 4F]  (4G) {${}^{4}\widetilde{G}$};
        \node [, right =1.0cm of 4G]  (dots) {$\cdots$};

          \draw[arrow]   (1G.east)  --(2F.west) node [pos=0.1,above] {}; 
          \draw[arrow]   (2F.east)  --(2G.west) node [pos=0.1,above] {}; 
          \draw[arrow]   (2G.east)  --(3F.west) node [pos=0.1,above] {};
          \draw[arrow]   (3F.east)  --(3G.west) node [pos=0.1,above] {};
          \draw[arrow]   (3G.east)  --(4F.west) node [pos=0.1,above] {};
          \draw[arrow]   (4F.east)  --(4G.west) node [pos=0.1,above] {};
          \draw[arrow]   (4G.east)--(dots.west);

    \end{tikzpicture}
    \caption{Flow chart for the recursive computation in (\ref{GM_wienerchaos}) and (\ref{FM_wienerchaos}).
        }\label{FlowChart}
\end{figure} 

\subsection{Elements of Malliavin Calculus}

Let $\mathcal{S}$ denote the set of all random variables $S$ of the form
$$s(W(f_{1}),\ldots,W(f_{n})),$$
where $n\geq1$, $s:\mathbb{R}^{n} \rightarrow \mathbb{R}$ is a $C^{\infty}$-function such that
$s$ and its partial derivatives have at most polynomial growth, and $f_{i}\in \overline{H}$, $i = 1,\ldots,n$.
We denote by $L^{2}(\Omega)$ the set of $\mathcal{F}_{W}$-measurable random variables
whose second moments exist.
The space $\mathcal{S}$ is dense in $L^{2}(\Omega)$ {\cite[Lemma 2.3.1]{nourdin2012normal}}.

\begin{Definition}
For every integers $n,p\geq1$, the $p$-th Malliavin derivative of
$$S=s(W(f_{1}),\ldots,W(f_{n}))\in \mathcal{S}$$
with respect to $W$ is defined by
\begin{align}\label{def:Dp}
D^{p}S = \overset{n}{\underset{i_{1},\ldots,i_{p}=1}{\sum}}\frac{\partial^{p}s}{\partial x_{i_{1}}\cdots \partial x_{i_{p}}}(W(f_{1}),\ldots,W(f_{n}))f_{i_{1}}\otimes \cdots \otimes f_{i_{p}}.
\end{align}
\end{Definition}

Because the sum in (\ref{def:Dp}) runs over all partial derivatives, $D^{p}S$ belongs to $L^{2}(\Omega,\overline{H}^{\odot p})$.
\begin{Definition}
For every integer $p\geq1$, we denote by Dom($\delta^{p}$) the subset of $L^{2}(\Omega,\overline{H}^{\otimes p})$
composed of those elements $u$ such that there exists a constant $c>0$ satisfying
\begin{equation}\label{def:dom_delta}
|\mathbb{E}\langle D^{p}S,u\rangle_{\overline{H}^{\otimes p}}|\leq c\sqrt{\mathbb{E}[S^{2}]}\ \textup{for all}\ S\in \mathcal{S}.
\end{equation}
\end{Definition}
For each $u\in \textup{Dom}(\delta^{p})$, (\ref{def:dom_delta}) and the Riesz representation theorem
imply that
there exists a unique element in $L^{2}(\Omega)$, denoted by $\delta^{p}(u)$,
such that
\begin{equation}\label{def:delta}
\mathbb{E}\left[ S\delta^{p}(u)\right] = \mathbb{E}\langle D^{p}S,u\rangle_{\overline{H}^{\otimes p}}\ \textup{for all}\ S\in \mathcal{S}.
\end{equation}
The operator $\delta^{p}: \textup{Dom}(\delta^{p})\subset L^{2}(\Omega,\overline{H}^{\otimes p})\rightarrow L^{2}(\Omega)$ is called the divergence operator of order $p$.
By \cite[Section 2.7]{nourdin2012normal}, we know that $\overline{H}^{\otimes p}\subset\ \textup{Dom}(\delta^{p})$
and $\delta^{p}(u)$ coincides with the $p$-fold
Wiener-It$\hat{\textup{o}}$ integrals of the function $u\in \overline{H}^{\otimes p}$ for any $p\geq1$;
that is,
$$\delta^{p}(u)=I_{p}(u)=\int_{\mathbb{R}^{p}}^{'}u(\lambda_{1},\ldots,\lambda_{p})W(d\lambda_{1})\cdots W(d\lambda_{p}).$$

The Malliavin derivative and divergence of multiple Wiener-It$\hat{\textup{o}}$ integrals have the following properties.
The proof can be {found in}
\cite[Proposition 2.7.4,Proposition 2.8.8]{nourdin2012normal}.
\begin{Lemma}[\cite{nourdin2012normal}]\label{lemma:deltaD} 
For every integer $p\geq1$ and $u\in \overline{H}^{\otimes p}$,
\begin{equation*}
\delta D I_{p}(u) = p I_{p}(u)
\end{equation*}
and
\begin{align}\label{lemma:DIp}
DI_{p}(u)
= pI_{p-1}(\tilde{u}),
\end{align}
where $\tilde{u}$ is the canonical symmetrization of $u$ defined through (\ref{def:canonical_symmetrization}) and
\begin{align*}\notag
I_{p-1}(\tilde{u})(\cdot) = \int^{'}_{\mathbb{R}^{p-1}}\tilde{u}(\lambda_{1},\ldots,\lambda_{p-1},\cdot) W(d\lambda_{1})\cdots W(d\lambda_{p-1}).
\end{align*}
On the other hand, for all $r\in \mathbb{N}$, $p\in \mathbb{N}$ and $q\in[1,\infty]$,
\begin{align}\label{norm_estimate1}
\|I_{p}(f)\|_{\mathbb{D}^{r,q}}\leq c_{r,p,q} \|f\|_{\overline{H}^{\otimes p}},\ f\in \overline{H}^{\odot p},
\end{align}
where $c_{r,p,q}>0$ is an universal constant and
\begin{align}\label{norm_estimate2}
\|I_{p}(f)\|_{\mathbb{D}^{r,q}}:=\left(\mathbb{E}[|I_{p}(f)|^{q}]
+\mathbb{E}[\|DI_{p}(f)\|_{\overline{H}}^{q}]+\cdots
+\mathbb{E}[\|D^{r}I_{p}(f)\|_{\overline{H}^{\otimes r}}^{q}
]\right)^{1/q}.
\end{align}

\end{Lemma}

We note that \eqref{norm_estimate1} is a special case of the Meyer's inequality \cite[Proposition 1.5.7]{nualart2006malliavin}.

\section{Main Results}\label{sec:mainresult}
As mentioned in Section \ref{sec:introduction}, we will present large-scale limit theorems about the non-Gaussian random process $U[j_{1:M}]X$
arising from STQNs of stationary Gaussian processes.

\subsection{Convergence Rate of STQNs with Random Inputs}
Let $\mathcal{N}$ be a standard normal random variable and
$\sigma_{j_{M}}^{2}=\mathbb{E}U[j_{1:M}]X(t).$
The total variation {(TV)} distance between the distributions of $\frac{1}{\sigma^{2}_{j_{M}}}U[j_{1:M}]X(2^{j_{M}}t)$ and the square of $\mathcal{N}$
is defined by
\begin{align*}
d_{\textup{TV}}\left(\frac{1}{\sigma_{j_{M}}^{2}}U[j_{1:M}]X(2^{j_{M}}t),\,\mathcal{N}^{2}\right)
=\underset{
h:\mathbb{R}\rightarrow[0,1]}{\textup{sup}}\left|\mathbb{E}\left[h\Big(\frac{1}{\sigma_{j_{M}}^{2}}U[j_{1:M}]X(2^{j_{M}}t)\Big)\right]-
\mathbb{E}\left[h(\mathcal{N}^{2})\right]\right|,
\end{align*}
where the supremum is taken over Borel functions $h$ taking values in $[0,1]$.
Because the composition function $h(|\cdot|^{2})$ also takes values in $[0,1]$,
\begin{align}\label{TV:U<TV:T}
d_{\textup{TV}}\left(\frac{1}{\sigma_{j_{M}}^{2}}U[j_{1:M}]X(2^{j_{M}}t),\,\mathcal{N}^{2}\right)
\leq
\underset{
h:\mathbb{R}\rightarrow[0,1]}{\textup{sup}}\left|\mathbb{E}\left[h\Big(\frac{1}{\sigma_{j_{M}}}T[j_{1:M}]X(2^{j_{M}}t)\Big)\right]-
\mathbb{E}\left[h(\mathcal{N})\right]\right|.
\end{align}
In the proposition below, we apply
Stein's method and the Malliavin calculus to the right hand side of (\ref{TV:U<TV:T})
to get an upper bound for $d_{\textup{TV}}\big(\frac{1}{\sigma_{j_{M}}^{2}}U[j_{1:M}]X(2^{j_{M}}t),\,\mathcal{N}^{2}\big)$.

\begin{Proposition}\label{prop:general_estimate}
For every stationary Gaussian process $X$, integer $M\geq2$, and $j_{1:M}=[j_{1},\ldots,j_{M}]\in \mathbb{Z}^{M}$,
the {TV} distance between the distributions of the non-Gaussian variable
$\frac{1}{\sigma^{2}_{j_{M}}}U[j_{1:M}]X(t)$ and the square of the standard normal random variable $\mathcal{N}$
has an upper bound
\begin{align}\label{prop:general_upper_bound}
d_{\textup{TV}}\left(\frac{1}{\sigma_{j_{M}}^{2}}U[j_{1:M}]X(2^{j_{M}}t),\,\mathcal{N}^{2}\right)\leq
\frac{2}{\sigma^{2}_{j_{M}}}\mathcal{U}_{j_{M}}(t,t),
\end{align}
where
\begin{align}\notag
&\mathcal{U}_{j_{M}}(s,t)=
 \overset{2^{M-2}-1}{\underset{r=0}{\sum}}(2^{M-1}-2r)\hspace{-0.2cm}
\overset{2^{M-1}-2r-1}{\underset{\ell=1}{\sum}}\hspace{-0.3cm}(\ell-1)!\binom{2^{M-1}-2r-1}{\ell-1}^{2}(2^{M}-4r-2\ell)!
\\\notag&{\times}\|{}^{M}\!F^{(2^{M-1}-2r)}_{2^{j_{M}}s}\otimes_{\ell}
{}^{M}\!F^{(2^{M-1}-2r)}_{2^{j_{M}}t}\|_{\overline{H}^{\otimes 2^{M}-4r-2\ell}}
\\\notag+&\overset{2^{M-2}-1}{\underset{\begin{subarray}{l}
r,r'=0\\
r\neq r'
\end{subarray}}{\sum}}(2^{M-1}-2r')\hspace{-0.2cm}
\overset{2^{M-1}-2(r\vee r')}{\underset{\ell=1}{\sum}}
\hspace{-0.3cm}(\ell-1)!\binom{2^{M-1}-2r-1}{\ell-1}\binom{2^{M-1}-2r'-1}{\ell-1}(2^{M}-2r-2r'-2\ell)!
\\\label{def:U_upperbound}&{\times}\|{}^{M}\!F^{(2^{M-1}-2r)}_{2^{j_{M}}s}\otimes_{\ell}
{}^{M}\!F^{(2^{M-1}-2r')}_{2^{j_{M}}t}\|_{\overline{H}^{\otimes 2^{M}-2r-2r'-2\ell}}
\end{align}
for $s,t\in \mathbb{R}$, and the family of functions $\{{}^{M}\!F^{(2^{M-1}-2r)}_{2^{j_{M}}t}\mid r=0,1,\ldots,2^{M-2}-1\}$ is defined through the recursive formula (\ref{GM_wienerchaos}) and (\ref{FM_wienerchaos}).
\end{Proposition}

In practice,
the scale parameters are set with $j_{m+1}\geq j_{m}$, where $m\in \mathbb{N}$, for the continuous wavelet transform in the scattering network
in order to use the neuron $U[j_{m+1}]$ in the $(m+1)$-th layer to extract larger-scale features
from the outputs of the neuron $U[j_{m+1}]$ in the $m$-th layer. {This is supported by the fact that the spectrum of the magnitude of the continuous wavelet transform is mainly supported on the lower frequency area. See \cite{mallat2012group} for more discussion. 
}
For fixed $j_{1},j_{2},\ldots,j_{M-1}\in \mathbb{Z}$,
the behavior of $\mathcal{U}_{j_{M}}$ and $\sigma^{2}_{j_{M}}$ when $j_{M}$ is sufficiently large
{will be} discussed in Proposition \ref{lemma:estimate_norm_fpqr}, whose proof relies on the following estimate{, and its proof is given in Section \ref{sec:estimate:G}}.
\begin{Lemma}\label{lemma:G_ineq}
For every $M\in \mathbb{N}$, there exists a constant $C_{M}>0$, which {depends on $j_{1},...,j_{M-1}$ but is} independent of $j_{M}$ and $t$, such that
\begin{equation}\label{G_ineq}
\left|{}^{M}\!\widetilde{G}^{(2^{M}-2\ell)}_{t}(\lambda_{1:2^{M}-2\ell})\right|\leq C_{M}
\overset{2^{M}-2\ell}{\underset{k=1}{\prod}}\sqrt{f_{X\star \psi_{j_{1}}}(\lambda_{k})}
\end{equation}
for $\ell = 0,1,\ldots,2^{M-1}-1$. 
\end{Lemma}

\begin{Proposition}\label{lemma:estimate_norm_fpqr}

(a) If the parameter $\alpha$ in Assumption \ref{Assumption:1:wavelet} and
the parameter $\beta$ in Assumption \ref{Assumption:2:spectral} satisfy $2\alpha+\beta\geq1$,
then for any $j_{1},\ldots,j_{M-1}\in \mathbb{Z}$, integer $M\geq2$, $r,r'\in\{0,1,\ldots,2^{M-2}-1\}$, and $s,t\in \mathbb{R}$,
\begin{align*}
2^{\frac{3}{2}j_{M}}\left\|{}^{M}\!F^{(2^{M-1}-2r)}_{2^{j_{M}}s}\otimes_{\ell}
{}^{M}\!F^{(2^{M-1}-2r')}_{2^{j_{M}}t}\right\|_{\overline{H}^{\otimes 2^{M}-2r-2r'-2\ell}}
\end{align*}
converges when $j_{M}\rightarrow\infty$, where $\ell\in\{1,2,\ldots,2^{M-1}-2r-1\}$ for the case $r=r'$ and
$\ell\in\{1,2,\ldots,2^{M-1}-2(r\vee r')\}$ for the case $r\neq r'$. It further implies that
$\underset{j_{M}\rightarrow\infty}{\lim}2^{\frac{3}{2}j_{M}}\mathcal{U}_{j_{M}}(s,t)$ exists and is finite{, where $\mathcal{U}_{j_{M}}$ is defined in \eqref{def:U_upperbound}}.

(b) Under the same condition, i.e., $2\alpha+\beta\geq1$,
\begin{align*}
\underset{j_{M}\rightarrow\infty}{\lim}2^{j_{M}}\sigma^{2}_{j_{M}}=
\overset{2^{M-2}-1}{\underset{r=0}{\sum}} M_{r}!\  c^{(M_{r})}\ \|\hat{\psi}\|^{2},
\end{align*}
where $M_{r} = 2^{M-1}-2r$ and
\begin{align}\label{def:c2m}
c^{(M_{r})}= \int_{\mathbb{R}^{M_{r}-1}}\Big|
{}^{M-1}\widetilde{G}^{(M_{r})}_{0}(u_{1:M_{r}-1},-u^{+}_{1:M_{r}-1})
\Big|^{2}
du_{1:M_{r}-1}.
\end{align}

\end{Proposition}

The condition $2\alpha+\beta\geq 1$ is used for {ensuring} $\|f_{X\star \psi_{j_{1}}}\|_{\infty}<\infty$,
which allows us to apply the Lebesgue dominated convergence theorem in the proof.  
{Note that this condition could be easily achieved for widely used wavelets, including
the Daubechies wavelets, the Morlet wavelet, the Mexican hat wavelet, etc, since $\alpha\geq\frac{1}{2}$ holds for those wavelets.}
The proofs of Lemma \ref{lemma:G_ineq} and Proposition \ref{lemma:estimate_norm_fpqr} can be found in
Section \ref{sec:estimate:G} and Section \ref{proof:lemma:estimate_norm_fpqr}, respectively.
By {Propositions} \ref{prop:general_estimate} and \ref{lemma:estimate_norm_fpqr},
we get the first main result as follows.

\begin{Theorem}\label{thm:conv_rate_4rd}
Let $\psi$ be a real-valued mother wavelet satisfying Assumption \ref{Assumption:1:wavelet}, and $X$ be a stationary Gaussian process satisfying Assumption \ref{Assumption:2:spectral}.
If the parameter $\alpha$ in Assumption \ref{Assumption:1:wavelet} and
the parameter $\beta$ in Assumption \ref{Assumption:2:spectral} satisfy $2\alpha+\beta\geq1$,
then for each fixed $t \in \mathbb{R}$, integer $M\geq2$,
and $j_{1},\ldots,j_{M-1}\in \mathbb{Z}$, when $j_{M}$ is sufficiently large,
\begin{equation*}
d_{\textup{TV}}\left(\frac{1}{\sigma^{2}_{j_{M}}}U[j_{1:M}]X(t),\,\mathcal{N}^{2}\right) = O(2^{-\frac{1}{2}j_{M}}).
\end{equation*}
\end{Theorem}
The theorem above implies that for each fixed $t\in \mathbb{R}$,
\begin{equation*}
\frac{1}{\sigma^{2}_{j_{M}}}U[j_{1:M}]X(t) \Rightarrow \mathcal{N}^{2}
\end{equation*}
in the distribution sense when $j_{M}\rightarrow\infty$, and hence together with Proposition \ref{lemma:estimate_norm_fpqr} we have
\begin{equation*}
2^{j_{M}}U[j_{1:M}]X(t) \Rightarrow \left[\overset{2^{M-2}-1}{\underset{r=0}{\sum}} M_{r}!\  c^{(M_{r})}\ \|\hat{\psi}\|^{2}\right]\mathcal{N}^{2}
\end{equation*}
when $j_{M}\rightarrow\infty$.

\subsection{Multivariate Normal Approximation of $2^{\frac{j_{M}}{2}}T[j_{1:M}]X(2^{j_{M}}t)$}
In this section, we further show that
$2^{j_{M}}U[j_{1:M}]X$ converges in the finite dimensional distribution sense {and} derive the covariance function of the limiting process.
Our approach is based on a result in \cite{nourdin2010invariance}.
First of all, by the Wiener-It$\hat{\textup{o}}$ decomposition \cite[Definition 2.2.4]{nourdin2012normal}, every {$F\in L^2(\Omega, \mathcal{F}_{W}, \mathcal P)$} 
admits a unique expansion of the type
$$
F = \mathbb{E}[F]+\overset{\infty}{\underset{p=1}{\sum}}\ \textup{Proj}(F\mid \mathcal{H}_{p}),
$$
where $\textup{Proj}(F\mid \mathcal{H}_{p})\in \mathcal{H}_{p}$ and the above series converges in $L^{2}(\Omega,\mathcal{F}_{W},\mathcal{P})$.
\begin{Definition}{}
For any $F\in L^{2}{(\Omega,\mathcal{F}_{W},\mathcal{P})}$, the pseudo-inverse of the infinitesimal generator $L$ of the Ornstein-Uhlenbeck semigroup, {denoted as $L^{-1}$,} is defined as
\begin{equation}\label{def:inverseL}
L^{-1}F = -\overset{\infty}{\underset{p=1}{\sum}}\ \frac{1}{p}\ \textup{Proj}(F\mid \mathcal{H}_{p}).
\end{equation}
\end{Definition}

The details about the semigroup of the Ornstein-Uhlenbeck semigroup and its infinitesimal generator
can be {found in} \cite[Section 2.8]{nourdin2012normal}.
{We also need the following lemma coming from the multidimensional Stein's method.}

\begin{Lemma}[\cite{nourdin2010invariance}]\label{cited_thm:multiGA}
For any integer $d\geq2$, let $\mathbf{F}=(F_{1},\ldots,F_{d})^\top$ be a $\mathcal{F}_{W}$-measurable random vector such that $\mathbb{E}[|F_{i}|^{4}]+\mathbb{E}[\|DF_{i}\|_{\overline{H}}^{4}]<\infty$
and $\mathbb{E}[F_{i}]=0$ for $i=1,\ldots,d$. Let $[\mathbf{C}(i,j)]_{1\leq i,j\leq d}$ be a symmetric non-negative definite matrix in $\mathbb{R}^{d\times d}$,
and let $\mathcal{N}_{\mathbf{C}}$ be a $d$-dimensional normal random vector with mean zero and covariance matrix $\mathbf{C}$. Then, for any twice differentiable function $h:\mathbb{R}^{d}\rightarrow \mathbb{R}$ with
$$
\|h^{''}\|_{\infty}:=\underset{1\leq i,j\leq d}{\sup}\ \underset{x\in \mathbb{R}^{d}}{\sup}\ \Big|\frac{\partial^{2}h}{\partial x_{i} \partial x_{j}}(x_{1},\ldots,x_{d})\Big|<\infty,
$$
{we have}
\begin{align}\label{multidim_stein}
|\mathbb{E}[h(\mathbf{F})]-\mathbb{E}[h(\mathcal{N}_{\mathbf{C}})]|\leq \frac{1}{2}\|h^{''}\|_{\infty} \sqrt{\overset{d}{\underset{i,j=1}{\sum}}
\mathbb{E}\left[\left(\mathbf{C}(i,j)-\Big\langle DF_{j},-DL^{-1}F_{i}\Big\rangle\right)^{2}\right]}.
\end{align}
\end{Lemma}

{With Lemma \ref{cited_thm:multiGA}, we proceed to study $2^{j_{M}}U[j_{1:M}]X$ in the finite dimensional distribution sense. We start from checking conditions needed for Lemma \ref{cited_thm:multiGA}.} For any $t_{1},t_{2},\ldots,t_{d}\in \mathbb{R}$, let
\begin{equation}\label{def:Fi}
F_{i} := 2^{\frac{j_{M}}{2}}T[j_{1:M}]X(2^{j_{M}}t_{i})=2^{\frac{j_{M}}{2}}
\left[\overset{2^{M-2}-1}{\underset{r=0}{\sum}}I_{2^{M-1}-2r}\left({}^{M}\!F^{(2^{M-1}-2r)}_{2^{j_{M}}t_{i}}
\right)
\right]
\end{equation}
for $i=1,2,\ldots,d.$
{ By (\ref{norm_estimate1}) and (\ref{norm_estimate2}), the condition $\mathbb{E}[|F_{i}|^{4}]+\mathbb{E}[\|DF_{i}\|_{\overline{H}}^{4}]<\infty$ holds.} %
On the other hand, because $0\leq r\leq 2^{M-2}-1$ in the summation (\ref{def:Fi}), $F_{i}$ is a linear combination of Wiener chaos of order greater than or equal to 2,
which implies that $\mathbb{E}[F_{i}]=0$.
Set {the covariance matrix of $(F_1,\ldots,F_d)^\top$ as $\mathbf{C}_{j_{M}}\in \mathbb{R}^{d\times d}$; that is,}
\begin{equation}\label{def:Fi_cov}
\mathbf{C}_{j_{M}}(i,j) := \mathbb{E}[F_{i}F_{j}],
\end{equation}
{where $1\leq i,j\leq d$,} which satisfies the requirement of symmetric and non-negative definite for the matrix $\mathbf{C}$ in Lemma \ref{cited_thm:multiGA}.
Moreover, {we claim that}
it can be expressed as
\begin{align}\label{lemma:multicov:expansion}
\mathbf{C}_{j_{M}}(i,j)  =
2^{j_{M}}\overset{2^{M-2}-1}{\underset{r=0}{\sum}}
M_{r}!\ {}^{M}\!F^{(M_{r})}_{2^{j_{M}}t_{i}}\otimes_{M_{r}} {}^{M}\!F^{(M_{r})}_{2^{j_{M}}t_{j}},
\end{align}
where $M_{r} = 2^{M-1}-2r$.
The proof of (\ref{lemma:multicov:expansion}) can be found in Section \ref{sec:proof:lemma:multicov:expansion}.
Because all conditions in Lemma \ref{cited_thm:multiGA} are satisfied,
the inequality (\ref{multidim_stein}) holds with $\mathbf{C}=\mathbf{C}_{j_{M}}$.
{Next, we deal with the right hand side of (\ref{multidim_stein}). To this end,
we} apply the Malliavin calculus, the equality (\ref{lemma:multicov:expansion}) and Proposition \ref{lemma:estimate_norm_fpqr}.
The result is summarized as follows, {and its proof is relegated to Section \ref{proof:thm:equiv_F}.}

\begin{Theorem}\label{thm:equiv_F}
Let $\psi$ be a real-valued mother wavelet satisfying Assumption \ref{Assumption:1:wavelet}, and $X$ be a stationary Gaussian process satisfying Assumption \ref{Assumption:2:spectral}.
Given integers $d,M\geq2$, $j_{1},\ldots,{j_M}\in \mathbb{Z}$, and $t_{1},t_{2},\ldots,t_{d}\in \mathbb{R}$,
define
$$
\mathbf{F} = \left(2^{\frac{j_{M}}{2}}T[j_{1:M}]X(2^{j_{M}}t_{1}),\ldots,2^{\frac{j_{M}}{2}}T[j_{1:M}]X(2^{j_{M}}t_{d})\right)^\top {\in \mathbb{R}^d}
$$
and let $\mathbf{C}_{j_{M}}$ be its covariance matrix.
If the parameter $\alpha$ in Assumption \ref{Assumption:1:wavelet} and
the parameter $\beta$ in Assumption \ref{Assumption:2:spectral} satisfy $2\alpha+\beta\geq1$, then for any twice differentiable function $h:\mathbb{R}^{d}\rightarrow \mathbb{R}$ with
$
\|h^{''}\|_{\infty}<\infty,
$
we have
\begin{align}\label{multidim_stein_j3}
\left|\mathbb{E}[h(\mathbf{F})]-\mathbb{E}[h(\mathcal{N}_{\mathbf{C}_{j_{M}}})]\right|\leq
\frac{1}{2}\|h^{''}\|_{\infty} \overset{d}{\underset{i,j=1}{\sum}}2^{j_{M}}\mathcal{U}_{j_{M}}(t_{j},t_{i})=
\|h^{''}\|_{\infty}O(2^{-\frac{1}{2}j_{M}})
\end{align}
when $j_{M}$ is sufficiently large, where $\mathcal{N}_{\mathbf{C}_{j_{M}}}$ is a $d$-dimensional normal random vector with mean zero and covariance matrix $\mathbf{C}_{j_{M}}$.
\end{Theorem}

By considering the special cases
$h(\mathbf{x}) = \cos(\langle\zeta,\mathbf{x}\rangle)$ and
$h(\mathbf{x}) = \sin(\langle\zeta,\mathbf{x}\rangle)$, where $\mathbf{x},\zeta\in \mathbb{R}^{d}$,
(\ref{multidim_stein_j3}) implies that the distance between the characteristic function of the random vector $\mathbf{F}$
and that of $\mathcal{N}_{\mathbf{C}_{j_{M}}}$ converges to zero as $j_{M}\rightarrow\infty$.
More precisely,
\begin{equation}\label{inq:characteristic_1}
\Big|\mathbb{E}\left[e^{i\langle \zeta,\,\mathbf{F}\rangle}\right]-\mathbb{E}\left[e^{i\langle \zeta,\,\mathcal{N}_{\mathbf{C}_{j_{M}}}\rangle}\right]\Big| \leq |\zeta|^{2}O(2^{-\frac{1}{2}j_{M}}).
\end{equation}
The following proposition shows that the covariance matrix $\mathbf{C}_{j_{M}}$ converges as $j_{M}\rightarrow\infty$.

\begin{Proposition}\label{prop:cov_matrix_lim}
Let the assumptions and notation of Theorem \ref{thm:equiv_F} prevail. Then, for any $i,j\in\{1,2,\ldots,d\}$,
we have
\begin{align}\label{def:Cinfty_3rd}
&\mathbf{C}_{\infty}(i,j):=\underset{j_{M}\rightarrow\infty}{\lim}\mathbf{C}_{j_{M}}(i,j)
=
\left[\overset{2^{M-2}-1}{\underset{r=0}{\sum}}M_{r}!\ c^{(M_{r})}\right]\int_{\mathbb{R}}e^{i(t_{i}-t_{j})z}\left|\hat{\psi}(z)\right|^{2}\ dz,
\end{align}
where $M_{r} = 2^{M-1}-2r$, the constant $c^{(M_{r})}$ is defined in (\ref{def:c2m}). 

\end{Proposition}

The proof of Proposition \ref{prop:cov_matrix_lim} is in Section \ref{proof:prop:cov_matrix_lim}.
Let $\mathcal{N}_{\mathbf{C}_{\infty}}$ be a $d$-dimensional normal random vector with mean zero and covariance matrix $\mathbf{C}_{\infty}$.
Because both $\mathcal{N}_{\mathbf{C}_{j_{M}}}$ and $\mathcal{N}_{\mathbf{C}_{\infty}}$ are normal random vectors with mean 0,
Proposition \ref{prop:cov_matrix_lim} implies that for any $\zeta=[\zeta_{1},\ldots,\zeta_{d}]\in \mathbb{R}^{d}$
\begin{equation}\label{char_normal_vector}
\mathbb{E}\left[e^{i\langle \zeta,\,\mathcal{N}_{\mathbf{C}_{j_{M}}}\rangle}\right]
\rightarrow
\mathbb{E}\left[e^{i\langle \zeta,\,\mathcal{N}_{\mathbf{C}_{\infty}}\rangle}\right]
\end{equation}
when $j_{M}\rightarrow\infty$.
By combining (\ref{inq:characteristic_1}) and (\ref{char_normal_vector}), we get
\begin{align*}\notag
&\Big|\mathbb{E}\left[e^{i\langle \zeta,\,\mathbf{F}\rangle}\right]-\mathbb{E}\left[e^{i\langle \zeta,\,\mathcal{N}_{\mathbf{C}_{\infty}}\rangle}\right]\Big|
\\\leq& \Big|\mathbb{E}\left[e^{i\langle \zeta,\,\mathbf{F}\rangle}\right]-\mathbb{E}\left[e^{i\langle \zeta,\,\mathcal{N}_{\mathbf{C}_{j_{M}}}\rangle}\right]\Big|+
\Big|\mathbb{E}\left[e^{i\langle \zeta,\,\mathcal{N}_{\mathbf{C}_{j_{M}}}\rangle}\right]-\mathbb{E}\left[e^{i\langle \zeta,\,\mathcal{N}_{\mathbf{C}_{\infty}}\rangle}\right]\Big|
\rightarrow0
\end{align*}
when $j_{M}\rightarrow\infty$.
Therefore, for any $i,j\in\{1,2,\ldots,d\}$, when $j_{M}\rightarrow\infty$,
\begin{align*}
\mathbb{E}\left[e^{i\langle \zeta,\,\mathbf{F}\rangle}\right]\rightarrow\mathbb{E}\left[e^{i\langle \zeta,\,\mathcal{N}_{\mathbf{C}_{\infty}}\rangle}\right]
\end{align*}
for all $\zeta\in \mathbb{R}^{d}$, where $\mathcal{N}_{\mathbf{C}_{\infty}}$ is a $d$-dimensional normal random vector with mean zero and covariance matrix $\mathbf{C}_{\infty}$ given {in} (\ref{def:Cinfty_3rd}).
It implies that
the rescaled random process
$\{2^{\frac{j_{M}}{2}}T[j_{1:M}]X(2^{j_{M}}t)\mid t\in \mathbb{R}\}$
converges to a Gaussian process, denoted by $\{G(t)\mid t\in \mathbb{R}\}$, in the finite dimensional distribution sense and
\begin{equation*}
\textup{Cov}\left(G(t_{1}),G(t_{2})\right)=\left[\overset{2^{M-2}-1}{\underset{r=0}{\sum}}M_{r}!\ c^{(M_{r})}\right]\int_{\mathbb{R}}e^{i(t_{1}-t_{2})z}\Big|\hat{\psi}(z)\Big|^{2}\ dz,\ \ t_{1},t_{2}\in \mathbb{R}.
\end{equation*}
By the continuous mapping theorem, we get a corresponding result for
the rescaled random process
$\{2^{j_{M}}U[j_{1:M}]X(2^{j_{M}}t)\mid t\in \mathbb{R}\}$
as follows.

\begin{Theorem}\label{thm:equiv_F_plus}
Let the assumptions and notation of Theorem \ref{thm:equiv_F} prevail.
The rescaled random process
$\{2^{j_{M}}U[j_{1:M}]X(2^{j_{M}}t)\mid t\in \mathbb{R}\}$
converges to the square of the stationary Gaussian process $\{G(t)\mid t\in \mathbb{R}\}$ in the finite dimensional distribution sense.
Moreover, the covariance function of the limiting process $G^{2}$ has the following expression
\begin{align}\notag
\textup{Cov}\left(G^{2}(t_{1}),G^{2}(t_{2})\right)
=&\, 2\left[\overset{2^{M-2}-1}{\underset{r=0}{\sum}}M_{r}!\ c^{(M_{r})}\right]^{2}\left[\int_{\mathbb{R}}e^{i(t_{1}-t_{2})z}|\hat{\psi}(z)|^{2}\ dz\right]^{2}
\\\label{cov:U}=&\,
2\left[\textup{Cov}\left(G(t_{1}),G(t_{2})\right)\right]^{2}
\end{align}
for $t_{1},t_{2}\in \mathbb{R},$ where $M_{r} = 2^{M-1}-2r$ and the constant $c^{(M_{r})}$ is defined in (\ref{def:c2m}).
\end{Theorem}
The derivation of (\ref{cov:U}) can be found in Section \ref{sec:proof:finalresult}.
Theorem \ref{thm:equiv_F_plus} implies that when $j_{M}\rightarrow\infty$,
$$2^{j_{M}}U[j_{1:M}]X(2^{j_{M}}t)\Rightarrow
\left[\overset{2^{M-2}-1}{\underset{r=0}{\sum}} M_{r}!\  c^{(M_{r})}\ \|\hat{\psi}\|^{2}\right]
\chi^{2}(t)
$$
in the finite dimensional distribution sense,
where $\{\chi^{2}(t)\mid t\in \mathbb{R}\}$ is a chi-square process with one degree of freedom
and
\begin{align}\notag
\textup{Cov}\left(\chi^{2}(t_{1}),\chi^{2}(t_{2})\right)
=& 2 \left[\int_{\mathbb{R}}e^{i(t_{1}-t_{2})z}\left(\frac{|\hat{\psi}(z)|}{\|\hat{\psi}\|_{2}}\right)^{2}\ dz\right]^{2}.
\end{align}

\begin{Remark}
It is worth mentioning that in our previous work \cite{liu2020central}, we considered
the second-order {ST} with stationary Gaussian processes and Gaussian processes with stationary increments
as inputs{, where we apply} the Wiener chaos expansion and the Feynman diagram method
to derive {the central and non-central limit theorem, and the delicate interaction between mother wavelets and activation functions are explored}.
Different from \cite{liu2020central}, 
{the analysis techniques in this work} allows us to get a general form for the Wiener-It$\hat{\textup{o}}$ decomposition
of the outputs of STQNs without depth limitation, {and the current approach is} non-asymptotic, i.e., the Stein's method is applied
to estimate the distance between the outputs of STQNs and their scaling limits.
For commonly used wavelet functions,
the current work provides information about the convergence speed of the rescaled random processes
arising from STQNs, which cannot be obtained from our previous work.
Furthermore,
the technique used in the current work can be applied to study the same issue for polynomial activated STs
by using more complicated product formula to obtain the Wiener-It$\hat{\textup{o}}$ decomposition
of the corresponding outputs{, which is not explored in this work.}

\end{Remark}

\section{Proof}\label{sec:proof}

\subsection{Proof of Lemma \ref{lemma:recursiveTU}}\label{sec:proof:lemma:recursiveTU}

In the following, we prove (\ref{TM_wienerchaos}) and (\ref{UM_wienerchaos}) by the mathematical induction method.

$\bullet$ When $M=2$, not only (\ref{GM_wienerchaos}) and (\ref{FM_wienerchaos}) can be simplified as (\ref{2G2F}) and (\ref{material_T2}) respectively,
but also (\ref{TM_wienerchaos}) and (\ref{UM_wienerchaos}) become
(\ref{T2_wienerchaos}) and (\ref{U2_wienerchaos}).
Hence, the statement in Lemma \ref{lemma:recursiveTU} holds for $M=2$.

$\bullet$ Suppose that the statement in Lemma \ref{lemma:recursiveTU} holds for $M=m$, where $m\geq2.$

$\bullet$ Now we start to prove the statement in Lemma \ref{lemma:recursiveTU} for $M=m+1$. By the definition (\ref{def:T}),
\begin{align}\notag
T[j_{1:m+1}]X(t) =&\, U[j_{1:m}]X\star\psi_{j_{m+1}}(t)
\\\label{induction1}=&\,\int_{\mathbb{R}}\overset{2^{m-1}}{\underset{\ell=0}{\sum}}I_{2^{m}-2\ell}\left({}^{m}\widetilde{G}^{(2^{m}-2\ell)}_{t^{'}}
\right)\psi_{j_{m+1}}(t-t^{'})dt',
\end{align}
where the second equality follows from the induction hypothesis.
By (\ref{time_relation_MG}),
\begin{align}\label{mG_psi}
\int_{\mathbb{R}}{}^{m}\widetilde{G}^{(2^{m}-2\ell)}_{t^{'}}(\lambda_{1:2^{m}-2\ell})\psi_{j_{m+1}}(t-t^{'})dt'
={}^{m}\widetilde{G}^{(2^{m}-2\ell)}_{t}(\lambda_{1:2^{m}-2\ell})\hat{\psi}_{j_{m+1}}(\lambda^{+}_{1:2^{m}-2\ell}).
\end{align}
We recognize that the right hand side of (\ref{mG_psi}) is just the right hand side of (\ref{FM_wienerchaos}), i.e,
${}^{m+1}F^{(2^{m}-2\ell)}_{t}(\lambda_{1:2^{m}-2\ell})$.
By the stochastic Fubini theorem \cite[Theorem 2.1]{pipiras2010regularization} with (\ref{induction1}) and (\ref{mG_psi}), we get
\begin{align}\label{induction2}
T[j_{1:m+1}]X(t) =
\overset{2^{m-1}}{\underset{\ell=0}{\sum}}I_{2^{m}-2\ell}\left({}^{m+1}F^{(2^{m}-2\ell)}_{t}
\right),
\end{align}
that is, (\ref{TM_wienerchaos}) with $M=m+1$. 
The validation of the application of the stochastic Fubini theorem to get (\ref{induction2}) can be found in \ref{appendix:fubini}.
Next, by the relationship (\ref{U=squareT}) and (\ref{induction2}),
\begin{align}\notag
&U[j_{1:m+1}]X(t) = \left[T[j_{1:m+1}]X(t)\right]^{2}
=\left[\overset{2^{m-1}}{\underset{r=0}{\sum}}I_{2^{m}-2r}\left({}^{m+1}F^{(2^{m}-2r)}_{t}
\right)\right]^{2}
\\\notag=&\,
\overset{2^{m-1}}{\underset{r,r'=0}{\sum}}
I_{2^{m}-2r}\left({}^{m+1}F^{(2^{m}-2r)}_{t}\right)
I_{2^{m}-2r'}\left({}^{m+1}F^{(2^{m}-2r')}_{t}\right)
\\\notag=&\,
\overset{2^{m-1}}{\underset{r,r'=0}{\sum}}
\overset{2^{m}-2(r\vee r')}{\underset{p=0}{\sum}}
p!\binom{2^{m}-2r}{p}\binom{2^{m}-2r'}{p}\ I_{2^{m+1}-2(r+r'+p)}\left({}^{m+1}F^{(2^{m}-2r)}_{t}
\otimes_{p} {}^{m+1}F^{(2^{m}-2r')}_{t}\right),
\end{align}
where the last equality follows from the product formula, i.e., Lemma \ref{lemma:itoformula}.
Denote $\ell = r+r'+p$. Because $0\leq p \leq 2^{m}-2(r\vee r')$,
we have $r+r'\leq \ell \leq 2^{m}-|r-r^{'}|$. Hence,
\begin{align}\notag
U[j_{1:m+1}]X(t) =&
\overset{2^{m-1}}{\underset{r,r'=0}{\sum}}\ \overset{2^{m}}{\underset{\ell=0}{\sum}}
1_{\{r+r'\leq \ell \leq 2^{m}-|r-r^{'}|\}}
(\ell-r-r')!\binom{2^{m}-2r}{\ell-r-r'}\binom{2^{m}-2r'}{\ell-r-r'}
\\\notag&\times I_{2^{m+1}-2\ell}\left({}^{m+1}F^{(2^{m}-2r)}_{t}
\otimes_{\ell-r-r'} {}^{m+1}F^{(2^{m}-2r')}_{t}\right),
\end{align}
where the indicator function
$1_{\{r+r'\leq \ell \leq 2^{m}-|r-r^{'}|\}}$
can be ignored by setting
\begin{equation}\notag
\binom{2^{m}-2r}{\ell-r-r'}\binom{2^{m}-2r'}{\ell-r-r'} = 0\ \textup{if}\ \ell-r-r'<0\
\textup{or}\ \min\left\{2^{m}-2r,2^{m}-2r'\right\}< \ell-r-r'.
\end{equation}
By the linearity of the Wiener integrals,
\begin{align}\notag
U[j_{1:m+1}]X(t) =&\,\overset{2^{m}}{\underset{\ell=0}{\sum}}I_{2^{m+1}-2\ell}\left(
\overset{2^{m-1}}{\underset{r,r'=0}{\sum}}
(\ell-r-r')!\binom{2^{m}-2r}{\ell-r-r'}\binom{2^{m}-2r'}{\ell-r-r'}\right.
\\\notag&\left.\times\ \  {}^{m+1}F^{(2^{m}-2r)}_{t}
\otimes_{\ell-r-r'} {}^{m+1}F^{(2^{m}-2r')}_{t}\right).
\end{align}
We recognize that the function inside the integrand of the $(2^{m+1}-2\ell)$-fold Wiener integrals
above coincides with the right hand side of (\ref{GM_wienerchaos}) with $M=m+1$,
i.e.,
${}^{m+1}\widetilde{G}^{(2^{m+1}-2\ell)}_{t}$.
Hence,
\begin{align}\notag
U[j_{1:m+1}]X(t) =&\,\overset{2^{m}}{\underset{\ell=0}{\sum}}I_{2^{m+1}-2\ell}\left(
{}^{m+1}\widetilde{G}^{(2^{m+1}-2\ell)}_{t}
\right),
\end{align}
which is just (\ref{UM_wienerchaos}) with $M=m+1$.
Therefore, we finish the inductive step and the proof.

\subsection{Proof of Proposition \ref{prop:general_estimate}}

Because the distribution of $T[j_{1:M}]X(2^{j_{M}}t)$ does not depend on the time variable,
we will neglect the time variable
of $T[j_{1:M}]X$ temporarily.
The Stein's equation associated with $h$ is the ordinary differential equation
\begin{equation}\label{stein_equation}
f^{'}(x)-xf(x) = h(x)-\mathbb{E}[h(\mathcal{N})].
\end{equation}
Let $f_{h}$ be the solution of (\ref{stein_equation}) with $\underset{|x|\rightarrow\infty}{\lim} f_{h}(x)e^{-\frac{x^{2}}{2}}=0$.
According to \cite[Proposition 3.2.2]{nourdin2012normal},
\begin{align*}
f_{h}(x)  = e^{\frac{x^{2}}{2}}\int_{-\infty}^{x} \left\{h(y)-\mathbb{E}\left[h(\mathcal{N})\right]\right\}e^{-\frac{y^{2}}{2}}dy.
\end{align*}
By (\ref{stein_equation}),
\begin{align}\notag
&\mathbb{E}\left[h\left(\frac{1}{\sigma_{j_{M}}}T[j_{1:M}]X\right)\right]-\mathbb{E}\left[h\left(\mathcal{N}\right)\right]
\\\label{E-E}=&\,
\mathbb{E}\left[f_{h}^{'}\left(\frac{1}{\sigma_{j_{M}}}T[j_{1:M}]X\right)\right]
-
\mathbb{E}\left[\frac{1}{\sigma_{j_{M}}}T[j_{1:M}]X\ f_{h}\left(\frac{1}{\sigma_{j_{M}}}T[j_{1:M}]X\right)\right].
\end{align}
By (\ref{TM_wienerchaos}), i.e.,
$T[j_{1:M}]X(2^{j_{M}}t)
=\overset{2^{M-2}-1}{\underset{r=0}{\sum}}I_{2^{M-1}-2r}\left({}^{M}\!F^{(2^{M-1}-2r)}_{2^{j_{M}}t}
\right)$, and Lemma \ref{lemma:deltaD}, we have
\begin{align*}
T[j_{1:M}]X(2^{j_{M}}t) = \delta D\left(\overset{2^{M-2}-1}{\underset{r=0}{\sum}}\frac{1}{2^{M-1}-2r}I_{2^{M-1}-2r}\left({}^{M}\!F^{(2^{M-1}-2r)}_{2^{j_{M}}t}\right)\right).
\end{align*}
Hence,
\begin{align}\notag
&\mathbb{E}\left[\frac{1}{\sigma_{j_{M}}}T[j_{1:M}]X\ f_{h}\left(\frac{1}{\sigma_{j_{M}}}T[j_{1:M}]X\right)\right]
\\\notag=&\,\frac{1}{\sigma_{j_{M}}}\mathbb{E}\left[\delta D\left(\overset{2^{M-2}-1}{\underset{r=0}{\sum}}\frac{1}{2^{M-1}-2r}I_{2^{M-1}-2r}\left({}^{M}\!F^{(2^{M-1}-2r)}_{2^{j_{M}}t}\right)\right)
f_{h}\left(\frac{1}{\sigma_{j_{M}}}T[j_{1:M}]X\right)\right]
\\\label{app:deltaD}=&\,\frac{1}{\sigma^{2}_{j_{M}}}\mathbb{E}\left[\Big\langle D\left(\overset{2^{M-2}-1}{\underset{r=0}{\sum}}\frac{1}{M_{r}}I_{M_{r}}\left({}^{M}\!F^{(M_{r})}_{2^{j_{M}}t}\right)\right),
f_{h}^{'}\left(\frac{1}{\sigma_{j_{M}}}T[j_{1:M}]X\right)DT[j_{1:M}]X\Big\rangle\right],
\end{align}
where the last equality follows from (\ref{def:delta}) and the chain rule.
By substituting (\ref{app:deltaD}) into (\ref{E-E}), we get
\begin{align}\notag
&\mathbb{E}\left[h\left(\frac{1}{\sigma_{j_{M}}}T[j_{1:M}]X\right)\right]-\mathbb{E}\left[h(\mathcal{N})\right]
\\\notag=&\,
\mathbb{E}\left\{f_{h}^{'}\left(\frac{1}{\sigma_{j_{M}}}T[j_{1:M}]X\right)
\left[1-\frac{1}{\sigma^{2}_{j_{M}}}\Big\langle D\left(\overset{2^{M-2}-1}{\underset{r=0}{\sum}}\frac{1}{M_{r}}I_{M_{r}}\left({}^{M}\!F^{(M_{r})}_{2^{j_{M}}t}\right)\right),\,DT[j_{1:M}]X
\Big\rangle
\right]\right\}.
\end{align}
By \cite[Theorem 3.3.1]{nourdin2012normal}, if the Borel function $h$ is assumed to take values in $[0,1]$,
the solution of Stein's equation associated with $h$ satisfies $\|f_{h}\|_{\infty}\leq \sqrt{\pi/2}$ and $\|f_{h}^{'}\|_{\infty}\leq 2$.
Hence,
\begin{align}\label{|E-E|}
&d_{\textup{TV}}\left(\frac{1}{\sigma_{j_{M}}}T[j_{1:M}]X,\,\mathcal{N}\right)
\leq
2\mathbb{E}\Big[
\Big|1-\frac{1}{\sigma^{2}_{j_{M}}}\Big\langle D\left(\overset{2^{M-2}-1}{\underset{r=0}{\sum}}\frac{1}{M_{r}}I_{M_{r}}\left({}^{M}\!F^{(M_{r})}_{2^{j_{M}}t}\right)\right),DT[j_{1:M}]X
\Big\rangle
\Big|\Big].
\end{align}
By applying (\ref{lemma:DIp}) in Lemma \ref{lemma:deltaD} to the right hand side of (\ref{|E-E|}),
\begin{align}\label{|E-E|v2}
d_{\textup{TV}}\left(\frac{1}{\sigma_{j_{M}}}T[j_{1:M}]X,\,\mathcal{N}\right)
\leq
2\mathbb{E}\left\{
\Big|1-\frac{1}{\sigma^{2}_{j_{M}}}
\mathcal{I}
\Big|\right\},
\end{align}
where
\begin{equation*}
\mathcal{I}:=\Big\langle \overset{2^{M-2}-1}{\underset{r=0}{\sum}}I_{M_{r}-1}\left({}^{M}\!F^{(M_{r})}_{2^{j_{M}}t}\right),
\overset{2^{M-2}-1}{\underset{r=0}{\sum}}M_{r}I_{M_{r}-1}\left({}^{M}\!F^{(M_{r})}_{2^{j_{M}}t}\right)\Big\rangle.
\end{equation*}
Note that
\begin{align}\notag
\mathbb{E}\left[\mathcal{I}\right]
=\overset{2^{M-2}-1}{\underset{r=0}{\sum}} M_{r}!\ \ \|{}^{M}\!F^{(M_{r})}_{t}\|_{\overline{H}^{\otimes M_{r}}}^{2}
= \sigma_{j_{M}}^{2}.
\end{align}
By applying Jensen's inequality to (\ref{|E-E|v2}),
we get
\begin{align}\label{TV_var}
d_{\textup{TV}}\left(\frac{1}{\sigma_{j_{M}}}T[j_{1:M}]X,\,\mathcal{N}\right)
\leq
2\left\{\mathbb{E}\left[
\Big|1-\frac{1}{\sigma^{2}_{j_{M}}}\mathcal{I}
\Big|^{2}\right]\right\}^{\frac{1}{2}}
=2\sqrt{\textup{Var}\left(\frac{1}{\sigma^{2}_{j_{M}}}\mathcal{I}\right)}.
\end{align}
Before computing the variance in (\ref{TV_var}), we expand $\mathcal{I}$ by Lemma \ref{lemma:itoformula} as follows
\begin{align}\notag
\mathcal{I}
=&\,\sigma_{j_{M}}^{2}+
\overset{2^{M-2}-1}{\underset{r=0}{\sum}}M_{r}
\overset{M_{r}-1}{\underset{\ell=1}{\sum}}(\ell-1)!\binom{M_{r}-1}{\ell-1}^{2}
I_{2M_{r}-2\ell}
\left({}^{M}\!F^{(M_{r})}_{2^{j_{M}}t}\otimes_{\ell}
{}^{M}\!F^{(M_{r})}_{t}\right)
\\\notag&+\overset{2^{M-2}-1}{\underset{\begin{subarray}{l}
r,r'=0\\
r\neq r'
\end{subarray}}{\sum}}M_{r'}
\overset{M_{r}\wedge M_{r'}}{\underset{\ell=1}{\sum}}(\ell-1)!\binom{M_{r}-1}{\ell-1}\binom{M_{r'}-1}{\ell-1}
I_{M_{r}+M_{r'}-2\ell}\left({}^{M}\!F^{(M_{r})}_{2^{j_{M}}t}\otimes_{\ell}
{}^{M}\!F^{(M_{r'})}_{2^{j_{M}}t}\right).
\end{align}
By the Minkowski inequality,
\begin{align}\notag
&\sqrt{\textup{Var}\left(\mathcal{I}\right)} = \mathbb{E}\left\{\left[\Big|\mathcal{I}-\sigma_{j_{M}}^{2}\Big|^{2}\right]\right\}^{\frac{1}{2}}
\\\label{sqrt_V_I}\leq& \overset{2^{M-2}-1}{\underset{r=0}{\sum}}M_{r}
\overset{M_{r}-1}{\underset{\ell=1}{\sum}}(\ell-1)!\binom{M_{r}-1}{\ell-1}^{2}
\left\{\mathbb{E}\left[\Big|I_{2M_{r}-2\ell}
\left({}^{M}\!F^{(M_{r})}_{2^{j_{M}}t}\otimes_{\ell}
{}^{M}\!F^{(M_{r})}_{2^{j_{M}}t}\right)\Big|^{2}\right]\right\}^{\frac{1}{2}}
\\\notag&+\overset{2^{M-2}-1}{\underset{\begin{subarray}{l}
r,r'=0\\
r\neq r'
\end{subarray}}{\sum}}M_{r'}
\overset{M_{r}\wedge M_{r'}}{\underset{\ell=1}{\sum}}(\ell-1)!\binom{M_{r}-1}{\ell-1}\binom{M_{r'}-1}{\ell-1}
\left\{\mathbb{E}\left[\Big|I_{M_{r}+M_{r'}-2\ell}\left({}^{M}\!F^{(M_{r})}_{2^{j_{M}}t}\otimes_{\ell}
{}^{M}\!F^{(M_{r'})}_{2^{j_{M}}t}\right)\Big|^{2}\right]\right\}^{\frac{1}{2}}.
\end{align}
For the expectations above, the isometry property of the multiple Wiener-It$\hat{\textup{o}}$ integrals implies that
\begin{align*}
\mathbb{E}\left[\Big|I_{M_{r}+M_{r'}-2\ell}\left({}^{M}\!F^{(M_{r})}_{2^{j_{M}}t}\otimes_{\ell}
{}^{M}\!F^{(M_{r'})}_{2^{j_{M}}t}\right)\Big|^{2}\right]
=
(M_{r}+M_{r'}-2\ell)!\ \  \left\|{}^{M}\!F^{(M_{r})}_{2^{j_{M}}t}\otimes_{\ell}
{}^{M}\!F^{(M_{r'})}_{2^{j_{M}}t}\right\|^{2}.
\end{align*}
Therefore, the upper bound in (\ref{prop:general_upper_bound})
follows from combining (\ref{TV_var}) and (\ref{sqrt_V_I}).

\subsection{Proof of Lemma \ref{lemma:G_ineq}}\label{sec:estimate:G}
We prove this lemma by mathematical induction.
For $M=1$, (\ref{def:1F2}) implies that (\ref{G_ineq}) holds with $C_{1}=1$.
Suppose that (\ref{G_ineq}) holds for $M=m-1$, where $m$ is an integer greater than or equal to 2.
By (\ref{GM_wienerchaos}),
\begin{align}\label{Gm_wienerchaos}
{}^{m}\widetilde{G}^{(2^{m}-2\ell)}_{t}
=\overset{2^{m-2}-1}{\underset{r,r^{'}=0}{\sum}}
(\ell-r-r')!\binom{2^{m-1}-2r}{\ell-r-r'}\binom{2^{m-1}-2r'}{\ell-r-r'}{}^{m}\!F^{(2^{m-1}-2r)}_{t}\widetilde{\otimes}_{\ell-r-r'} {}^{m}\!F^{(2^{m-1}-2r')}_{t}.
\end{align}
From (\ref{FM_wienerchaos}), we know that for $r = 0,1,\ldots,2^{m-2}-1$,
\begin{align}\notag
\Big|{}^{m}\!F^{(2^{m-1}-2r)}_{t}(x_{1},\ldots,x_{2^{m-1}-2r})\Big|=&\, \Big|{}^{m-1}\widetilde{G}^{(2^{m-1}-2r)}_{t}(x_{1},\ldots,x_{2^{m-1}-2r})\hat{\psi}_{j_{m}}(x_{1},\ldots,x_{2^{m-1}-2r})\Big|
\\\notag\leq&\,
\|\hat{\psi}\|_{\infty}\Big|{}^{m-1}\widetilde{G}^{(2^{m-1}-2r)}_{t}(x_{1},\ldots,x_{2^{m-1}-2r})\Big|
\\\notag\leq&\,
\|\hat{\psi}\|_{\infty}\left[C_{m-1}
\overset{2^{m-1}-2r}{\underset{k=1}{\prod}}\sqrt{f_{X\star \psi_{j_{1}}}(x_{k})} \right],
\end{align}
where the last inequality is obtained from the assumption, i.e., (\ref{Gm_wienerchaos}) with $M=m-1$.
Hence,
\begin{align}\notag
\Big|{}^{m}\!F^{(2^{m-1}-2r)}_{t}\widetilde{\otimes}_{\ell-r-r'} {}^{m}\!F^{(2^{m-1}-2r')}_{t}(\lambda_{1:2^{m}-2\ell})\Big|
\leq
\|\hat{\psi}\|^{2}_{\infty}C_{m-1}^{2}\|f_{X\star \psi_{j_{1}}}\|_{1}^{\ell-r-r'}
\overset{2^{m}-2\ell}{\underset{k=1}{\prod}}\sqrt{f_{X\star \psi_{j_{1}}}(\lambda_{k})}.
\end{align}
By substituting the inequality above into (\ref{Gm_wienerchaos}),
we get
\begin{align}\notag
\Big|{}^{m}\widetilde{G}^{(2^{m}-2\ell)}_{t}(\lambda_{1:2^{m}-2\ell})\Big|
\leq&\,\left[\overset{2^{m-2}-1}{\underset{r,r^{'}=0}{\sum}}
(\ell-r-r')!\binom{2^{m-1}-2r}{\ell-r-r'}\binom{2^{m-1}-2r'}{\ell-r-r'}\|f_{X\star \psi_{j_{1}}}\|_{1}^{\ell-r-r'}\right]
\\\notag&\times\left[\|\hat{\psi}\|^{2}_{\infty}C_{m-1}^{2}
\overset{2^{m}-2\ell}{\underset{k=1}{\prod}}\sqrt{f_{X\star \psi_{j_{1}}}(\lambda_{k})}\right].
\end{align}
By the principle of mathematical induction, (\ref{G_ineq}) holds for all $M\in \mathbb{N}$.

\subsection{Proof of Proposition \ref{lemma:estimate_norm_fpqr}}\label{proof:lemma:estimate_norm_fpqr}

{\bf (a)}
For $r,r'\in\{0,1,\ldots,2^{M-2}-1\}$,
by the recursive relationship (\ref{FM_wienerchaos}),
\begin{align}\notag
&\left\|{}^{M}\!F^{(M_{r})}_{2^{j_{M}}s}\otimes_{\ell}{}^{M}\!F^{(M_{r'})}_{2^{j_{M}}t}\right\|^{2}_{\overline{H}^{\otimes M_{r}+M_{r'}-2\ell}}
\\\notag=&\int_{\mathbb{R}^{M_{r}+M_{r'}-2\ell}}
\Big[\int_{\mathbb{R}^{\ell}} {}^{M-1}\widetilde{G}^{(M_{r})}_{2^{j_{M}}s}(\tau_{1:\ell},\lambda_{1:M_{r}-\ell})\hat{\psi}_{j_{M}}(\tau^{+}_{1:\ell}+\lambda^{+}_{1:M_{r}-\ell})
\\\notag&\times{}^{M-1}\widetilde{G}^{(M_{r'})}_{2^{j_{M}}t}(-\tau_{1:\ell},\lambda_{M_{r}-\ell+1:M_{r}+M_{r'}-2\ell})
\hat{\psi}_{j_{M}}(-\tau^{+}_{1:\ell}+\lambda^{+}_{M_{r}-\ell+1:M_{r}+M_{r'}-2\ell}){d\tau_{1:\ell}}\Big]
\\\notag&\times\Big[\overline{\int_{\mathbb{R}^{\ell}}{}^{M-1}\widetilde{G}^{(M_{r})}_{2^{j_{M}}s}(\eta_{1:\ell},\lambda_{1:M_{r}-\ell})\hat{\psi}_{j_{M}}(\eta^{+}_{1:\ell}+\lambda^{+}_{1:M_{r}-\ell})}
\\\notag&\overline{\times{}^{M-1}\widetilde{G}^{(M_{r'})}_{2^{j_{M}}t}(-\eta_{1:\ell},\lambda_{M_{r}-\ell+1:M_{r}+M_{r'}-2\ell})
\hat{\psi}_{j_{M}}(-\eta^{+}_{1:\ell}+\lambda^{+}_{M_{r}-\ell+1:M_{r}+M_{r'}-2\ell}){d\eta_{1:\ell}}} \Big]
d\lambda_{1:M_{r}+M_{r'}-2\ell},
\end{align}
where $\ell\in\{1,\ldots,M_{r}-1\}$ {when} $r=r'$, and $\ell\in\{1,\ldots,M_{r}\wedge M_{r'}\}$ {when} $r\neq r'$.
By considering the change of variables
\begin{align}\notag
\left\{\begin{array}{l}
u_{i}=\tau_{i},\ i=1,\ldots,\ell-1\ \ \textup{for}\ \ell>1,
\\
v_{i}=\eta_{i},\ i=1,\ldots,\ell-1\ \ \textup{for}\ \ell>1,
\\
w_{i} = \lambda_{i},\ i =1,\ldots,M_{r}+M_{r'}-2\ell-1,
\\
x=2^{j_{M}}\left(\tau^{+}_{1:\ell}+\lambda^{+}_{1:M_{r}-\ell}\right),
\\
y=2^{j_{M}}\left(-\tau^{+}_{1:\ell}+\lambda^{+}_{M_{r}-\ell+1:M_{r}+M_{r'}-2\ell}\right),
\\
z=2^{j_{M}}\left(\eta^{+}_{1:\ell}+\lambda^{+}_{1:M_{r}-\ell}\right)
\end{array}
\right.
\end{align}
and noting that $\hat{\psi}_{j_{M}}(\cdot) = \hat{\psi}(2^{j_{M}}\cdot)$,
we have
\begin{align}\notag
&2^{3j_{M}}\left\|{}^{M}\!F^{(M_{r})}_{2^{j_{M}}s}\otimes_{\ell}{}^{M}\!F^{(M_{r'})}_{2^{j_{M}}t}\right\|^{2}_{\overline{H}^{\otimes M_{r}+M_{r'}-2\ell}}
\\\notag=&\,\int_{\mathbb{R}^{M_{r}+M_{r'}}}
{}^{M-1}\widetilde{G}^{(M_{r})}_{2^{j_{M}}s}(u_{1:\ell-1},u^{*}_{\ell},w_{1:M_{r}-\ell})\nonumber\\
\notag&
{\times}{}^{M-1}\widetilde{G}^{(M_{r'})}_{2^{j_{M}}t}(-u_{1:\ell-1},-u^{*}_{\ell},w_{M_{r}-\ell+1:M_{r}+M_{r'}-2\ell-1},w^{*}_{M_{r}+M_{r'}-2\ell})
\\\notag&\times\overline{{}^{M-1}\widetilde{G}^{(M_{r})}_{2^{j_{M}}s}(v_{1:\ell-1},v_{\ell}^{*},w_{1:M_{r}-\ell})}\\
\notag &{\times}
\overline{{}^{M-1}\widetilde{G}^{(M_{r'})}_{2^{j_{M}}t}(-v_{1:\ell-1},-v^{*}_{\ell},w_{M_{r}-\ell+1:M_{r}+M_{r'}-2\ell-1},w^{*}_{M_{r}+M_{r'}-2\ell})
}
\\\label{eq:check_limit}&\times
\hat{\psi}(x)\hat{\psi}(y)\overline{\hat{\psi}(z)\hat{\psi}(x+y-z)}\ du_{1:\ell-1}\ dv_{1:\ell-1}\ dw_{1:M_{r}+M_{r'}-2\ell-1}
\ dx\ dy\ dz,
\end{align}
where
$u^{*}_{\ell} = 2^{-j_{M}}x-u_{1}-\cdots-u_{\ell-1}-w_{1}-\cdots-w_{M_{r}-\ell}$,
$v^{*}_{\ell} = 2^{-j_{M}}z-v_{1}-\cdots-v_{\ell-1}-w_{1}-\cdots-w_{M_{r}-\ell}$, and
$w^{*}_{M_{r}+M_{r'}-2\ell} = 2^{-j_{M}}x+2^{-j_{M}}y-w_{1}-\cdots-w_{M_{r}+M_{r'}-2\ell-1}$.
In the following, we show that
the limit of the integrand on the right hand side of the equality (\ref{eq:check_limit}) exists when $j_{M}\rightarrow\infty$.
\begin{itemize}
\item For the second component of the integrand in (\ref{eq:check_limit}), {by (\ref{time_relation_MG}),} 
\begin{align*}
&{}^{M-1}\widetilde{G}^{(M_{r'})}_{2^{j_{M}}t}(-u_{1:\ell-1},-u^{*}_{\ell},w_{M_{r}-\ell+1:M_{r}+M_{r'}-2\ell-1},w^{*}_{M_{r}+M_{r'}-2\ell})
\\=&\,\textup{exp}\left\{i2^{j_{M}}t\left[-u^{+}_{1:\ell-1}-u^{*}_{\ell}+w^{+}_{M_{r}-\ell+1:M_{r}+M_{r'}-2\ell-1}+w^{*}_{M_{r}+M_{r'}-2\ell}
\right]\right\}
\\&\times{}^{M-1}\widetilde{G}^{(M_{r'})}_{0}(-u_{1:\ell-1},-u^{*}_{\ell},w_{M_{r}-\ell+1:M_{r}+M_{r'}-2\ell-1},w^{*}_{M_{r}+M_{r'}-2\ell})
\\=&\,e^{ity}
\times{}^{M-1}\widetilde{G}^{(M_{r'})}_{0}(-u_{1:\ell-1},-u^{*}_{\ell},w_{M_{r}-\ell+1:M_{r}+M_{r'}-2\ell-1},w^{*}_{M_{r}+M_{r'}-2\ell})
\\\rightarrow&\,
e^{ity}
{\times} {}^{M-1}\widetilde{G}^{(M_{r'})}_{0}(-u_{1:\ell-1},u^{+}_{1:\ell-1}+w^{+}_{1:M_{r}-\ell},
w_{M_{r}-\ell+1:M_{r}+M_{r'}-2\ell-1},-w^{+}_{1:M_{r}+M_{r'}-2\ell-1})
\end{align*}
{when $j_M\to \infty$. 
}

\item For the fourth component of the integrand in (\ref{eq:check_limit}),
\begin{align*}
&{}^{M-1}\widetilde{G}^{(M_{r'})}_{2^{j_{M}}t}(-v_{1:\ell-1},-v^{*}_{\ell},w_{M_{r}-\ell+1:M_{r}+M_{r'}-2\ell-1},w^{*}_{M_{r}+M_{r'}-2\ell})
\\=&\,\textup{exp}\left\{i2^{j_{M}}t\left[-v^{+}_{1:\ell-1}-v^{*}_{\ell}+w^{+}_{M_{r}-\ell+1:M_{r}+M_{r'}-2\ell-1}+
w^{*}_{M_{r}+M_{r'}-2\ell}\right]\right\}
\\&\,\times{}^{M-1}\widetilde{G}^{(M_{r'})}_{0}(-v_{1:\ell-1},-v^{*}_{\ell},w_{M_{r}-\ell+1:M_{r}+M_{r'}-2\ell-1},w^{*}_{M_{r}+M_{r'}-2\ell})
\\=&\,
e^{it(x+y-z)}
{\times} {}^{M-1}\widetilde{G}^{(M_{r'})}_{0}(-v_{1:\ell-1},-v^{*}_{\ell},w_{M_{r}-\ell+1:M_{r}+M_{r'}-2\ell-1},w^{*}_{M_{r}+M_{r'}-2\ell})
\\\rightarrow&\,e^{it(x+y-z)}
{\times} {}^{M-1}\widetilde{G}^{(M_{r'})}_{0}(-v_{1:\ell-1},v^{+}_{1:\ell-1}+w^{+}_{1:M_{r}-\ell},
w_{M_{r}-\ell+1:M_{r}+M_{r'}-2\ell-1},-w^{+}_{1:M_{r}+M_{r'}-2\ell-1})
\end{align*}
{when $j_{M}\rightarrow\infty$.}

\item Similarly, for the first and third components of the integrand in (\ref{eq:check_limit}), we have
$$
{}^{M-1}\widetilde{G}^{(M_{r})}_{2^{j_{M}}s}(u_{1:\ell-1},u^{*}_{\ell},w_{1:M_{r}-\ell})\rightarrow
e^{isx}{\times} {}^{M-1}\widetilde{G}^{(M_{r})}_{0}(u_{1:\ell-1},-u^{+}_{1:\ell-1}-w^{+}_{1:M_{r}-\ell},w_{1:M_{r}-\ell})
$$
and
$$
{}^{M-1}\widetilde{G}^{(M_{r})}_{2^{j_{M}}s}(v_{1:\ell-1},v_{\ell}^{*},w_{1:M_{r}-\ell})
\rightarrow e^{isz} {\times} {}^{M-1}\widetilde{G}^{(M_{r})}_{0}(v_{1:\ell-1},-v^{+}_{1:\ell-1}-w^{+}_{1:M_{r}-\ell},w_{1:M_{r}-\ell})
$$
respectively when $j_{M}\rightarrow\infty$.
\end{itemize}

On the other hand,
by Lemma \ref{lemma:G_ineq}, {for the integrand in \eqref{eq:check_limit},} there exists a constant $C>0$, which is independent of $j_{M}$, such that
\begin{align}\notag
&\underset{x,y,z\in \mathbb{R}}{\sup}\Big|{}^{M-1}\widetilde{G}^{(M_{r})}_{2^{j_{M}}s}(u_{1:\ell-1},u^{*}_{\ell},w_{1:M_{r}-\ell}){\times}
{}^{M-1}\widetilde{G}^{(M_{r'})}_{2^{j_{M}}t}(-u_{1:\ell-1},-u^{*}_{\ell},w_{M_{r}-\ell+1:M_{r}+M_{r'}-2\ell-1},w^{*}_{M_{r}+M_{r'}-2\ell})
\\\notag&\times\overline{{}^{M-1}\widetilde{G}^{(M_{r})}_{2^{j_{M}}s}(v_{1:\ell-1},v_{\ell}^{*},w_{1:M_{r}-\ell})}{\times}
\overline{{}^{M-1}\widetilde{G}^{(M_{r'})}_{2^{j_{M}}t}(-v_{1:\ell-1},-v^{*}_{\ell},w_{M_{r}-\ell+1:M_{r}+M_{r'}-2\ell-1},w^{*}_{M_{r}+M_{r'}-2\ell})
}\Big|
\\\notag\leq&\,
C
\left[\underset{x,y,z\in \mathbb{R}}{\sup}f_{X\star\psi_{j_{1}}}(u_{\ell}^{*})f_{X\star\psi_{j_{1}}}(v_{\ell}^{*})f_{X\star\psi_{j_{1}}}(w_{M_{r}+M_{r'}-2\ell}^{*})\right]
\\\label{m<n_allcase}&\times\left[\overset{\ell-1}{\underset{k=1}{\prod}}f_{X\star\psi_{j_{1}}}(u_{k})f_{X\star\psi_{j_{1}}}(v_{k})\right]
\left[\overset{M_{r}+M_{r'}-2\ell-1}{\underset{k=1}{\prod}}f_{X\star\psi_{j_{1}}}(w_{k})\right].
\end{align}
If $2\alpha+\beta\geq1$, then $f_{X\star\psi_{j_{1}}}(\cdot) = 2^{2\alpha j_{1}}|C_{\hat{\psi}}(2^{j_{1}}\cdot)|^{2}C_{X}(\cdot)|\cdot|^{2\alpha+\beta-1}$ is not only integrable, but also bounded,
which implies that the right hand side of (\ref{m<n_allcase}) is integrable with respect to $(u_{1:\ell-1},
v_{1:\ell-1},w_{1:M_{r}+M_{r'}-2\ell-1})$.
By the Lebesgue dominated convergence theorem,
we know that {$\underset{j_{M}\rightarrow\infty}{\lim}2^{3j_{M}}\big\|{}^{M}F^{(M_{r})}_{2^{j_{M}}s}\otimes_{\ell}{}^{M}F^{(M_{r'})}_{2^{j_{M}}t}\big\|^{2}_{\overline{H}^{\otimes M_{r}+M_{r'}-2\ell}}$ exists and satisfies}
\begin{align}\notag
&\underset{j_{M}\rightarrow\infty}{\lim}2^{3j_{M}}\left\|{}^{M}F^{(M_{r})}_{2^{j_{M}}s}\otimes_{\ell}{}^{M}F^{(M_{r'})}_{2^{j_{M}}t}\right\|^{2}_{\overline{H}^{\otimes M_{r}+M_{r'}-2\ell}}
\\\notag\leq &\,C\|f_{X\star \psi_{j_{1}}}\|_{1}^{M_{r}+M_{r'}-3}
\|f_{X\star \psi_{j_{1}}}\|_{\infty}^{3}
\|\hat{\psi}\|_{\infty}^{3}
\int_{\mathbb{R}^{3}}\Big|\hat{\psi}(x)
\hat{\psi}(y)
\overline{\hat{\psi}(z)}
\overline{\hat{\psi}(x+y-z)}\Big|\
dx\ dy\ dz.
\end{align}
The proof of the statement (a) is completed.

{\bf (b)} For $r\in\{0,1,\ldots,2^{M-2}-1\}$, by (\ref{TM_wienerchaos}) and the orthogonal property of multiple Wiener-It$\hat{\textup{o}}$ integrals,
\begin{align}\label{proof:variance_limit_A}
\sigma_{j_{M}}^{2} = \mathbb{E}\Big|T[j_{1:M}]X(2^{j_{M}}t)\Big|^{2}
=\mathbb{E}\Big|\overset{2^{M-2}-1}{\underset{r=0}{\sum}}I_{M_{r}}\left({}^{M}\!F^{(M_{r})}_{2^{j_{M}}t}\right)\Big|^{2}
=\overset{2^{M-2}-1}{\underset{r=0}{\sum}} M_{r}! \|{}^{M}\!F^{(M_{r})}_{2^{j_{M}}t}\|_{\overline{H}^{\otimes M_{r}}}^{2}.
\end{align}
{From (\ref{time_relation_MF}),} we know that the value of $\|{}^{M}\!F^{(M_{r})}_{2^{j_{M}}t}\|_{\overline{H}^{\otimes M_{r}}}^{2}$
does not change along with the time variable.  Hence, we choose $t=0$.
From (\ref{FM_wienerchaos}), for each $r\in\{0,1,\ldots,2^{M-2}-1\}$,
\begin{align*}
\|{}^{M}\!F^{(M_{r})}_{0}\|_{\overline{H}^{\otimes M_{r}}}^{2}
=\int_{\mathbb{R}^{M_{r}}}
\Big|{}^{M-1}\widetilde{G}^{(M_{r})}_{0}(\lambda_{1:M_{r}})
\hat{\psi}_{j_{M}}(\lambda^{+}_{1:M_{r}})\Big|^{2}\ d\lambda_{1:M_{r}}.
\end{align*}
By changing of variables $u_{i} = \lambda_{i}$ for $i=1,\ldots,M_{r}-1$ and $z=2^{j_{M}}(\lambda_{1}+\cdots+\lambda_{M_{r}})$,
\begin{align*}
2^{j_{M}}\|{{}^{M}}\!F^{(M_{r})}_{0}\|_{\overline{H}^{\otimes M_{r}}}^{2}
= \int_{\mathbb{R}^{M_{r}}}
\Big|
{}^{M-1}\widetilde{G}^{(M_{r})}_{0}(u_{1:M_{r}-1},2^{-j_{M}}z-u^{+}_{1:M_{r}-1})\Big|^{2}
|\hat{\psi}(z)|^{2}
\ du_{1:M_{r}-1}\ dz.
\end{align*}
For the integrand of the integral above, by Lemma \ref{lemma:G_ineq},
there exists a constant $C_{M}$, which is independent of $j_{M}$, such that
\begin{align*}
\Big|
{}^{M-1}\widetilde{G}^{(M_{r})}_{0}(u_{1:M_{r}-1},2^{-j_{M}}z-u^{+}_{1:M_{r}-1})\Big|^{2}
\leq C^{2}_{M}f_{X\star\psi_{j_{1}}}(2^{-j_{M}}z-u^{+}_{1:M_{r}-1}) \overset{M_{r}-1}{\underset{k=1}{\prod}}f_{X\star\psi_{j_{1}}}(u_{k}).
\end{align*}
Under the assumption $2\alpha+\beta\geq1$, $f_{X\star\psi_{j_{1}}}$ is not only integrable, but also bounded.
Hence, we can {again} apply the Lebesgue dominated convergence theorem to get
\begin{align}\label{proof:variance_limit_B}
\underset{j_{M}\rightarrow\infty}{\lim}2^{j_{M}}\|{{}^{M}}\!F^{(M_{r})}_{0}\|_{H^{\otimes M_{r}}}^{2}
= c^{(M_{r})}\|\hat{\psi}\|^{2},
\end{align}
where $c^{(M_{r})}$ is defined in (\ref{def:c2m}).
The proof is finished by combining (\ref{proof:variance_limit_A}) and (\ref{proof:variance_limit_B}).

\subsection{Proof of the equality (\ref{lemma:multicov:expansion})}\label{sec:proof:lemma:multicov:expansion}
From the definition of $\mathbf{C}_{j_{M}}(i,j)$ in (\ref{def:Fi_cov}), i.e., $\mathbf{C}_{j_{M}}(i,j)  =  \mathbb{E}[F_{i}F_{j}]$
and the definition of $F_{i}$ in (\ref{def:Fi}), i.e.,
$F_{i} = 2^{\frac{j_{M}}{2}}T[j_{1:M}]X(2^{j_{M}}t_{i})$,
we have
\begin{align*}
\mathbf{C}_{j_{M}}(i,j)  = 2^{j_{M}}\mathbb{E}\Big[T[j_{1:M}]X(2^{j_{M}}t_{i})T[j_{1:M}]X(2^{j_{M}}t_{j})\Big].
\end{align*}
By replacing $T[j_{1:M}]X$ with its Wiener-It$\hat{\textup{o}}$ decomposition in (\ref{TM_wienerchaos}),
\begin{align}\notag
\mathbf{C}_{j_{M}}(i,j)  =& \,2^{j_{M}}\mathbb{E}\left\{
\left[\overset{2^{M-2}-1}{\underset{r=0}{\sum}}I_{M_{r}}\left({}^{M}\!F^{(M_{r})}_{2^{j_{M}}t_{i}}\right)
\right]
\left[\overset{2^{M-2}-1}{\underset{r'=0}{\sum}}I_{M_{r'}}\left({}^{M}\!F^{(M_{r'})}_{2^{j_{M}}t_{j}}\right)
\right]\right\}
\\\label{proof:cov_eq}=&\,
2^{j_{M}}\overset{2^{M-2}-1}{\underset{r=0}{\sum}}\mathbb{E}\left[
I_{M_{r}}\left({}^{M}\!F^{(M_{r})}_{2^{j_{M}}t_{i}}\right)
I_{M_{r}}\left({}^{M}\!F^{(M_{r})}_{2^{j_{M}}t_{j}}\right)
\right],
\end{align}
where the last equality follows from the orthogonal property of the multiple Wiener integrals.
For $r\in\{0,1,\ldots,2^{M-2}-1\}$,
by Lemma \ref{lemma:itoformula},
\begin{align}\notag
I_{M_{r}}\left({}^{M}\!F^{(M_{r})}_{2^{j_{M}}t_{i}}\right)
I_{M_{r}}\left({}^{M}\!F^{(M_{r})}_{2^{j_{M}}t_{j}}\right)
= \overset{M_{r}}{\underset{\ell=0}{\sum}}\ell!\binom{M_{r}}{\ell}\binom{M_{r}}{\ell}
I_{2M_{r}-2\ell}\left({}^{M}\!F^{(M_{r})}_{2^{j_{M}}t_{i}}\otimes_{\ell} {}^{M}\!F^{(M_{r})}_{2^{j_{M}}t_{j}}\right).
\end{align}
Because the expectation of the Wiener integrals is equal to zero except the constant term,
\begin{align}\label{proof:cov_eq_II}
\mathbb{E}\left[I_{M_{r}}\left({}^{M}\!F^{(M_{r})}_{2^{j_{M}}t_{i}}\right)
I_{M_{r}}\left({}^{M}\!F^{(M_{r})}_{2^{j_{M}}t_{j}}\right)\right]
=M_{r}!\ \ {}^{M}\!F^{(M_{r})}_{2^{j_{M}}t_{i}}\otimes_{M_{r}} {}^{M}\!F^{(M_{r})}_{2^{j_{M}}t_{j}}.
\end{align}
The equality (\ref{lemma:multicov:expansion}) follows from substituting (\ref{proof:cov_eq_II}) into (\ref{proof:cov_eq}).

\subsection{Proof of Theorem \ref{thm:equiv_F}}\label{proof:thm:equiv_F}
As mentioned before the statement of Theorem \ref{thm:equiv_F},
all conditions in Lemma \ref{cited_thm:multiGA} are satisfied
when the random vector $\mathbf{F}$ and the matrix $\mathbf{C}$ are given by (\ref{def:Fi}) and (\ref{def:Fi_cov}).
Hence, it suffices to study the behavior of
$\mathbb{E}\left[\left(\mathbf{C}_{j_{M}}(i,j)-\Big\langle DF_{j},-DL^{-1}F_{i}\Big\rangle\right)^{2}\right]$.

By (\ref{def:inverseL}) and (\ref{lemma:DIp}), we have
\begin{align}\notag
DF_{j} =\overset{2^{M-2}-1}{\underset{r=0}{\sum}}M_{r}I_{M_{r}-1}\left({}^{M}\!F^{(M_{r})}_{2^{j_{M}}t_{j}}
\right)
\end{align}
and
\begin{align}\notag
-DL^{-1}F_{i} = 2^{\frac{j_{M}}{2}}D\left[\overset{2^{M-2}-1}{\underset{r=0}{\sum}}\frac{1}{M_{r}}I_{M_{r}}\left({}^{M}\!F^{(M_{r})}_{2^{j_{M}}t_{i}}
\right)\right]
= 2^{\frac{j_{M}}{2}}\left[\overset{2^{M-2}-1}{\underset{r=0}{\sum}}I_{M_{r}-1}\left({}^{M}\!F^{(M_{r})}_{2^{j_{M}}t_{i}}
\right)\right].
\end{align}
Therefore,
\begin{align}\label{multidim_stein_inner}
\Big\langle DF_{j},-DL^{-1}F_{i}\Big\rangle = 2^{j_{M}}\Big\langle \overset{2^{M-2}-1}{\underset{r=0}{\sum}}M_{r}I_{M_{r}-1}\left({}^{M}\!F^{(M_{r})}_{2^{j_{M}}t_{j}}
\right),
\overset{2^{M-2}-1}{\underset{r'=0}{\sum}}I_{M_{r'}-1}\left({}^{M}\!F^{(M_{r'})}_{2^{j_{M}}t_{i}}
\right)\Big\rangle.
\end{align}
By Lemma \ref{lemma:itoformula},
(\ref{multidim_stein_inner}) can be rewritten as follows
\begin{align}\label{DFDLinverse}
\Big\langle DF_{j},-DL^{-1}F_{i}\Big\rangle =
2^{j_{M}}\overset{2^{M-2}-1}{\underset{r=0}{\sum}}M_{r}!\left[{}^{M}\!F^{(M_{r})}_{2^{j_{M}}t_{j}}\otimes_{M_{r}}{}^{M}\!F^{(M_{r})}_{2^{j_{M}}t_{i}}\right]
+\mathcal{I}_{1}+\mathcal{I}_{2},
\end{align}
where
\begin{align}\notag
\mathcal{I}_{1}= 2^{j_{M}}
\overset{2^{M-2}-1}{\underset{r=0}{\sum}}M_{r}\overset{M_{r}-1}{\underset{\ell=1}{\sum}}
(\ell-1)!\binom{M_{r}-1}{\ell-1}^{2}
 I_{2M_{r}-2\ell}\left({}^{M}\!F^{(M_{r})}_{2^{j_{M}}t_{j}}\otimes_{\ell}{}^{M}\!F^{(M_{r})}_{2^{j_{M}}t_{i}}
\right)
\end{align}
and
\begin{align}\notag
\mathcal{I}_{2}= 2^{j_{M}}
\overset{2^{M-2}-1}{\underset{\begin{subarray}{l}
r,r'=0\\
r\neq r'
\end{subarray}}{\sum}}M_{r}\overset{M_{r}\wedge M_{r'}}{\underset{\ell=1}{\sum}}(\ell-1)!\binom{M_{r}-1}{\ell-1}\binom{M_{r'}-1}{\ell-1}
 I_{M_{r}+M_{r'}-2\ell}\left({}^{M}\!F^{(M_{r})}_{2^{j_{M}}t_{j}}\otimes_{\ell}{}^{M}\!F^{(M_{r'})}_{2^{j_{M}}t_{i}}
\right).
\end{align}
We recognize that the first term in (\ref{DFDLinverse}) coincides with the expression (\ref{lemma:multicov:expansion}) of
$\mathbf{C}_{j_{M}}(i,j)$.
By the Minkowski inequality,
\begin{align}\notag
&\sqrt{\mathbb{E}\left[\left(\mathbf{C}_{j_{M}}(i,j)-\Big\langle DF_{j},-DL^{-1}F_{i}\Big\rangle\right)^{2}\right]}
= \sqrt{\mathbb{E}\left[|\mathcal{I}_{1}+ \mathcal{I}_{2}|^{2}\right]}
\\\notag
\leq&\,\sqrt{\mathbb{E}\left[|\mathcal{I}_{1}|^{2}\right]}+\sqrt{\mathbb{E}\left[|\mathcal{I}_{2}|^{2}\right]}
= 2^{j_{M}}\mathcal{U}_{j_{M}}(t_{j},t_{i}),
\end{align}
where $\mathcal{U}_{j_{M}}(\cdot,\cdot)$ is defined in (\ref{def:U_upperbound}).
By Proposition \ref{lemma:estimate_norm_fpqr}, we know that
$\underset{j_{M}\rightarrow\infty}{\lim}2^{\frac{3}{2}j_{M}}\mathcal{U}_{j_{M}}(\cdot,\cdot)$ exists.
Hence, for $i,j\in\{1,2,\ldots,d\}$, we get
\begin{align*}
\sqrt{\mathbb{E}\left[\left(\mathbf{C}_{j_{M}}(i,j)-\Big\langle DF_{j},-DL^{-1}F_{i}\Big\rangle\right)^{2}\right]}=O(2^{-\frac{1}{2}j_{M}})
\end{align*}
when $j_{M}$ is sufficiently large.


\subsection{Proof of Proposition \ref{prop:cov_matrix_lim}}\label{proof:prop:cov_matrix_lim}

By the equality (\ref{lemma:multicov:expansion}), it suffices to prove that
$
2^{j_{M}}\ {}^{M}\!F^{(M_{r})}_{2^{j_{M}}t_{i}}\otimes_{M_{r}} {}^{M}\!F^{(M_{r})}_{2^{j_{M}}t_{j}}
$
converges when $j_{M}\rightarrow\infty$ for each $r\in\{0,1,\ldots,2^{M-2}-1\}$.
By (\ref{FM_wienerchaos}),
\begin{align}\notag
&
{}^{M}\!F^{(M_{r})}_{2^{j_{M}}t_{i}}\otimes_{M_{r}} {}^{M}\!F^{(M_{r})}_{2^{j_{M}}t_{j}}
=\int_{\mathbb{R}^{M_{r}}}{}^{M}\!F^{(M_{r})}_{2^{j_{M}}t_{i}}(\lambda_{1:M_{r}})
\ \overline{{}^{M}\!F^{(M_{r})}_{2^{j_{M}}t_{j}}(\lambda_{1:M_{r}})}\ d\lambda_{1:M_{r}}
\\\notag
=&\int_{\mathbb{R}^{M_{r}}}
\left[{}^{M-1}\widetilde{G}^{(M_{r})}_{2^{j_{M}}t_{i}}(\lambda_{1:M_{r}})\hat{\psi}_{j_{M}}(\lambda^{+}_{1:M_{r}})\right]
\ \left[\ \overline{{}^{M-1}\widetilde{G}^{(M_{r})}_{2^{j_{M}}t_{j}}(\lambda_{1:M_{r}})\hat{\psi}_{j_{M}}(\lambda^{+}_{1:M_{r}})}\ \right]\ d\lambda_{1:M_{r}}
\\\notag=&
\int_{\mathbb{R}^{M_{r}}}
\left[{}^{M-1}\widetilde{G}^{(M_{r})}_{2^{j_{M}}t_{i}}(\lambda_{1:M_{r}})\ \overline{{}^{M-1}\widetilde{G}^{(M_{r})}_{2^{j_{M}}t_{j}}(\lambda_{1:M_{r}})}\ \right]
 \Big|\hat{\psi}_{j_{M}}(\lambda^{+}_{1:M_{r}})\Big|^{2}\ d\lambda_{1:M_{r}}.
\end{align}
From the recursive formula (\ref{GM_wienerchaos}) and (\ref{FM_wienerchaos}) with the initial term (\ref{def:1F2}),
we know that ${}^{M-1}\widetilde{G}^{(M_{r})}_{t}(\lambda_{1:M_{r}})$ satisfies
$$
{}^{M-1}\widetilde{G}^{(M_{r})}_{t}(\lambda_{1:M_{r}}) = {}^{M-1}\widetilde{G}^{(M_{r})}_{0}(\lambda_{1:M_{r}})\exp\left(it\lambda^{+}_{1:M_{r}}\right).
$$
Hence,
\begin{align}\notag
{}^{M}\!F^{(M_{r})}_{2^{j_{M}}t_{i}}\otimes_{M_{r}} {}^{M}\!F^{(M_{r})}_{2^{j_{M}}t_{j}}
=
\int_{\mathbb{R}^{M_{r}}}\exp\left(i2^{j_{M}}(t_{i}-t_{j})\lambda^{+}_{1:M_{r}}\right)
\Big|{}^{M-1}\widetilde{G}^{(M_{r})}_{0}(\lambda_{1:M_{r}})\Big|^{2}
 \Big|\hat{\psi}_{j_{M}}(\lambda^{+}_{1:M_{r}})\Big|^{2}\ d\lambda_{1:M_{r}}.
\end{align}
By changing of variables
\begin{align*}
\left\{\begin{array}{l}
u_{k} = \lambda_{k}\ \textup{for}\ k=1,2,\ldots,M_{r}-1,\\
z= 2^{j_{M}}\lambda^{+}_{1:M_{r}},
\end{array}\right.
\end{align*}
we get
\begin{align}\notag
&2^{j_{M}}\ {}^{M}\!F^{(M_{r})}_{2^{j_{M}}t_{i}}\otimes_{M_{r}} {}^{M}\!F^{(M_{r})}_{2^{j_{M}}t_{j}}
\\\notag=&
\int_{\mathbb{R}^{M_{r}}}e^{i(t_{i}-t_{j})z}
\Big|{}^{M-1}\widetilde{G}^{(M_{r})}_{0}(u_{1:M_{r}-1},2^{-j_{M}}z-u^{+}_{1:M_{r}-1})\Big|^{2}
 \Big|\hat{\psi}(z)\Big|^{2}\ du_{1:M_{r}-1}\ dz.
\end{align}
Under the condition $2\alpha+\beta\geq 1$, we apply Lemma \ref{lemma:G_ineq} and the Lebesgue dominated convergence theorem to get
\begin{align}\notag
&\underset{j_{M}\rightarrow\infty}{\lim}2^{j_{M}}{}^{M}\!F^{(M_{r})}_{2^{j_{M}}t_{i}}\otimes_{M_{r}} {}^{M}\!F^{(M_{r})}_{2^{j_{M}}t_{j}}
\\\notag=&
\int_{\mathbb{R}}e^{i(t_{i}-t_{j})z}\left\{\int_{\mathbb{R}^{M_{r}-1}}
\Big|{}^{M-1}\widetilde{G}^{(M_{r})}_{0}(u_{1:M_{r}-1},-u^{+}_{1:M_{r}-1})\Big|^{2}
 \ du_{1:M_{r}-1}\right\}\Big|\hat{\psi}(z)\Big|^{2}\ dz.
\end{align}
Note that the integral inside the curly brackets above
is just the constant $c^{(M_{r})}$ defined in (\ref{def:c2m}).

\subsection{Proof of Theorem \ref{thm:equiv_F_plus} {(derivation of (\ref{cov:U}))}}\label{sec:proof:finalresult}

It suffices to derive the covariance function of the limiting process $G^{2}$.
From Proposition \ref{prop:cov_matrix_lim},
the covariance function of the process $G$ has the spectral representation
\begin{equation*}
\textup{Cov}\left(G(t_{1}),G(t_{2})\right)=\kappa\int_{\mathbb{R}}e^{i(t_{1}-t_{2})z}\Big|\hat{\psi}(z)\Big|^{2}\ dz,
\end{equation*}
where $\kappa$ is a nonnegative constant defined as
$$
\kappa = \overset{2^{M-2}-1}{\underset{r=0}{\sum}}M_{r}!\ c^{(M_{r})}.
$$
By introducing
a complex-valued Gaussian random measure $W^{'}$ on $\mathbb{R}$, which
also satisfies (\ref{ortho}),
the Gaussian process $G$ can be expressed as a Wiener integral
\begin{equation*}
G(t)=\sqrt{\kappa}\int_{\mathbb{R}}e^{itz}\hat{\psi}(z)\ W^{'}(dz),\ t\in \mathbb{R}.
\end{equation*}
By Lemma \ref{lemma:itoformula},
\begin{align*}
G^{2}(t)=&\,\kappa\left(\int_{\mathbb{R}}e^{itz}\hat{\psi}(z)\ W^{'}(dz)\right)^{2}
\\=&\,\kappa \int_{\mathbb{R}}|\hat{\psi}(z)|^{2}dz +\kappa \int^{'}_{\mathbb{R}^{2}}e^{it(z_{1}+z_{2})}\hat{\psi}(z_{1})\hat{\psi}(z_{2})W^{'}(dz_{1})W^{'}(dz_{2}),
\end{align*}
that is,
\begin{align*}
G^{2}(t)-\mathbb{E}[G^{2}(t)]=\kappa \int^{'}_{\mathbb{R}^{2}}e^{it(z_{1}+z_{2})}\hat{\psi}(z_{1})\hat{\psi}(z_{2})W^{'}(dz_{1})W^{'}(dz_{2}).
\end{align*}
By applying Lemma \ref{lemma:itoformula} to the product $\left(G^{2}(t_{1})-\mathbb{E}[G^{2}(t_{1})]\right)\left(G^{2}(t_{2})-\mathbb{E}[G^{2}(t_{2})]\right)$
and noting that the expectation of Wiener integrals is equal to zero,
we obtain
\begin{align*}
\textup{Cov}\left(G^{2}(t_{1}),G^{2}(t_{2})\right)
=&\,2\kappa^{2} \int_{\mathbb{R}^{2}}e^{i(t_{1}-t_{2})(z_{1}+z_{2})}|\hat{\psi}(z_{1})|^{2}|\hat{\psi}(z_{2})|^{2}\ dz_{1}\ dz_{2}.
\\=&\,
2\left[\kappa \int_{\mathbb{R}}e^{i(t_{1}-t_{2})z}|\hat{\psi}(z)|^{2}\ dz\right]^{2}
\\=&\,
2\left[\textup{Cov}\left(G(t_{1}),G(t_{2})\right)\right]^{2}.
\end{align*}

\section*{Reference}
\bibliographystyle{elsarticle-num}
\bibliography{reference_STQN}

\numberwithin{equation}{section}
\begin{appendix}
\section{A validation of the application of the stochastic Fubini theorem to get (\ref{induction2})}\label{appendix:fubini}
First of all, we apply the mathematical induction method to prove that for each fixed $t\in \mathbb{R}$,
\begin{align}\label{G_L2_claim}
{}^{M}\widetilde{G}^{(2^{M}-2\ell)}_{t}\in L^{2}(\mathbb{R}^{2^{M}-2\ell})\ \textup{for}\ \ell = 0,1,...,2^{M-1}-1,
\end{align}
(${}^{M}\widetilde{G}^{(0)}_{t}$ is a constant)
and
\begin{equation}\label{F_L2_claim}
{}^{M}F^{(2^{M-1}-2r)}_{t}\in L^{2}(\mathbb{R}^{2^{M-1}-2r})\ \textup{for}\ r = 0,1,\ldots,2^{M-2}-1.
\end{equation}

$\bullet$ From the initial term (\ref{def:1F2}), we know that for each fixed $t\in \mathbb{R}$,
${}^{1}\widetilde{G}^{(2)}_{t}\in L^{2}(\mathbb{R}^{2})$
because the spectral density function $f_{X}$ is integrable and $\psi\in L^{1}(\mathbb{R})$ implies that $\hat{\psi}$ is a bounded function.
The initial term (\ref{def:1F2}) also shows that ${}^{1}\widetilde{G}^{(0)}_{t}$ is a constant.
Hence, (\ref{G_L2_claim}) holds when $M=1$.
On the other hand, when $M=2$, (\ref{FM_wienerchaos}) can be simplified as
\begin{equation*}
{}^{2}F^{(2)}_{t}(\lambda_{1:2})
=
{}^{1}\widetilde{G}^{(2)}_{t}(\lambda_{1:2})
\hat{\psi}_{j_{M}}(\lambda^{+}_{1:2}).
\end{equation*}
By the fact $\|\hat{\psi}\|_{\infty}<\infty$, we obtain
${}^{2}F^{(2)}_{t}\in L^{2}(\mathbb{R}^{2})$.
Hence, (\ref{F_L2_claim}) holds when $M=2$.
At the same time, (\ref{GM_wienerchaos}) can be simplified as
\begin{align*}
{}^{2}\widetilde{G}^{(4-2\ell)}_{t}
=
\ell!\binom{2}{\ell}\binom{2}{\ell}
\times {}^{2}F^{(2)}_{t}\widetilde{\otimes}_{\ell} {}^{2}F^{(2)}_{t}\ \textup{for}\ \ell=0,1,2.
\end{align*}
Clearly, for each fixed $t\in \mathbb{R}$, ${}^{2}\widetilde{G}^{(0)}_{t}$ is a constant.
For $\ell =0, 1$, by the Minkowski inequality and the Cauchy-Schwarz inequality,
\begin{align*}
\|{}^{2}F^{(2)}_{t}\widetilde{\otimes}_{\ell} {}^{2}F^{(2)}_{t}\|_{2}\leq \|{}^{2}F^{(2)}_{t}\|_{2}^{2}.
\end{align*}
It implies that ${}^{2}\widetilde{G}^{(4-2\ell)}_{t}\in L^{2}(\mathbb{R}^{4-2\ell})$
for $\ell = 0,1$.
In other words, the claim (\ref{G_L2_claim}) holds when $M=2$.

$\bullet$ Now, we suppose that (\ref{G_L2_claim}) and (\ref{F_L2_claim}) hold when $M=m$, where  $m\geq2.$

$\bullet$ When $M=m+1$, from (\ref{FM_wienerchaos}),  i.e.,
\begin{equation*}
{}^{m+1}F^{(2^{m}-2r)}_{t}(\lambda_{1:2^{m}-2r})
=
{}^{m}\widetilde{G}^{(2^{m}-2r)}_{t}(\lambda_{1:2^{m}-2r})
\hat{\psi}_{j_{m+1}}(\lambda^{+}_{1:2^{m}-2r})\ \textup{for}\ r=0,1,\ldots,2^{m-1}-1,
\end{equation*}
we have
\begin{equation*}
\|{}^{m+1}F^{(2^{m}-2r)}_{t}\|_{2}
\leq
\|{}^{m}\widetilde{G}^{(2^{m}-2r)}_{t}\|_{2}
\|\hat{\psi}_{j_{m+1}}\|_{\infty}\ \textup{for}\ r=0,1,\ldots,2^{m-1}-1,
\end{equation*}
By $\|\hat{\psi}\|_{\infty}<\infty$ and the induction hypothesis
${}^{m}\widetilde{G}^{(2^{m}-2r)}_{t}\in L^{2}(\mathbb{R}^{2^{m}-2r})$ for $r=0,1,\ldots,2^{m-1}-1$,
we know that
\begin{equation}\label{mF_L2}
{}^{m+1}F^{(2^{m}-2r)}_{t}\in L^{2}(\mathbb{R}^{2^{m}-2r})\ \textup{for}\ r=0,1,\ldots,2^{m-1}-1.
\end{equation}
On the other hand, by (\ref{GM_wienerchaos}),
\begin{align*}
{}^{m+1}\widetilde{G}^{(2^{m+1}-2\ell)}_{t}
= &\overset{2^{m-1}-1}{\underset{r,r^{'}=0}{\sum}}
(\ell-r-r')!\binom{2^{m}-2r}{\ell-r-r'}\binom{2^{m}-2r'}{\ell-r-r'}\nonumber\\
&\times {}^{m+1}F^{(2^{m}-2r)}_{t}\widetilde{\otimes}_{\ell-r-r'} {}^{m+1}F^{(2^{m}-2r')}_{t}
\end{align*}
for $\ell = 0,1,\ldots,2^{m}$.
It shows that ${}^{m+1}\widetilde{G}^{(0)}_{t}$ is a constant.
For $\ell = 0,1,\ldots,2^{m}-1$, by the Minkowski inequality and the Cauchy-Schwarz inequality,
\begin{align*}
\|{}^{m+1}\widetilde{G}^{(2^{m+1}-2\ell)}_{t}\|_{2}
\leq &\overset{2^{m-1}-1}{\underset{r,r^{'}=0}{\sum}}
(\ell-r-r')!\binom{2^{m}-2r}{\ell-r-r'}\binom{2^{m}-2r'}{\ell-r-r'}\nonumber\\
&\times \|{}^{m+1}F^{(2^{m}-2r)}_{t}\widetilde{\otimes}_{\ell-r-r'} {}^{m+1}F^{(2^{m}-2r')}_{t}\|_{2}
\\\leq &\overset{2^{m-1}-1}{\underset{r,r^{'}=0}{\sum}}
(\ell-r-r')!\binom{2^{m}-2r}{\ell-r-r'}\binom{2^{m}-2r'}{\ell-r-r'}\nonumber\\
&\times \|{}^{m+1}F^{(2^{m}-2r)}_{t}\|_{2}\ \|{}^{m+1}F^{(2^{m}-2r')}_{t}\|_{2}<\infty,
\end{align*}
where the last inequality follows from (\ref{mF_L2}).
Therefore, we finish the inductive step and the proof of the claim (\ref{G_L2_claim}) and (\ref{F_L2_claim}).

For $M\geq 1$, from (\ref{def:T}), (\ref{UM_wienerchaos}) and the mean-zero property of the wavelet,
\begin{align*}
&T[j_{1:M+1}]X(t) = \int_{\mathbb{R}} U[j_{1:M}]X(s) \psi_{j_{M+1}}(t-s)ds
\\=&\overset{2^{M-1}}{\underset{\ell=0}{\sum}}\int_{\mathbb{R}} I_{2^{M}-2\ell}\left({}^{M}\widetilde{G}^{(2^{M}-2\ell)}_{s}\right) \psi_{j_{M+1}}(t-s)ds
\\=&\overset{2^{M-1}-1}{\underset{\ell=0}{\sum}}\int_{\mathbb{R}} I_{2^{M}-2\ell}\left({}^{M}\widetilde{G}^{(2^{M}-2\ell)}_{s}\right) \psi^{+}_{j_{M+1}}(t-s)ds
\\&-\overset{2^{M-1}-1}{\underset{\ell=0}{\sum}}\int_{\mathbb{R}} I_{2^{M}-2\ell}\left({}^{M}\widetilde{G}^{(2^{M}-2\ell)}_{s}\right) \psi^{-}_{j_{M+1}}(t-s)ds,
\end{align*}
where $\psi_{j_{M+1}}^{+}$ (resp. $\psi_{j_{M+1}}^{-}$) is the nonnegative part (resp. the negative part) of $\psi_{j_{M+1}}$.
Before applying the stochastic Fubini theorem \cite[Theorem 2.1]{pipiras2010regularization} to change the order of the deterministic integral
$\int_{\mathbb{R}}\ \cdots\ \psi^{+}_{j_{M+1}}(t-s)ds$ and the Wiener integral $I_{2^{M}-2\ell}(\cdots)$, we need to verify the condition
\begin{equation*}
\int_{\mathbb{R}} \|{}^{M}\widetilde{G}^{(2^{M}-2\ell)}_{s}\|_{2}\  \psi^{+}_{j_{M+1}}(t-s)ds<\infty
\end{equation*}
This condition holds because (\ref{time_relation_MG}) implies that
$\|{}^{M}\widetilde{G}^{(2^{M}-2\ell)}_{s}\|_{2}$ does not depend on $s$ and $\psi^{+}_{j_{M+1}}\in L^{1}$.
Similarly,
$\int_{\mathbb{R}} \|{}^{M}\widetilde{G}^{(2^{M}-2\ell)}_{s}\|_{2}\  \psi^{-}_{j_{M+1}}(t-s)ds<\infty.$
Therefore, the stochastic Fubini theorem implies that
\begin{align*}
T[j_{1:M+1}]X(t)
=&\overset{2^{M-1}-1}{\underset{\ell=0}{\sum}} I_{2^{M}-2\ell}\left(\int_{\mathbb{R}}{}^{M}\widetilde{G}^{(2^{M}-2\ell)}_{s}\psi^{+}_{j_{M+1}}(t-s)ds -\int_{\mathbb{R}}{}^{M}\widetilde{G}^{(2^{M}-2\ell)}_{s}\psi^{-}_{j_{M+1}}(t-s)ds\right)
\\=&
\overset{2^{M-1}-1}{\underset{\ell=0}{\sum}} I_{2^{M}-2\ell}\left(
{}^{M+1}F^{(2^{M}-2\ell)}_{t}\right),
\end{align*}
where the last equality follows from the definition (\ref{FM_wienerchaos}).

\end{appendix}

\end{document}